# Optimal control of Goursat-Volterra systems.

S. A. Belbas
Mathematics Department
University of Alabama
Tuscaloosa, AL 35487-0350. USA.

e-mail: SBELBAS@BAMA.UA.EDU

Abstract. We analyze an optimal control problem for systems of integral equations of Volterra type with two independent variables. These systems generalize both, the hyperbolic control problems for systems of Goursat-Darboux type, and the optimal control of ordinary (i.e. with one independent variable) Volterra integral equations. We prove extremal principles akin to Pontryagin's maximum principle.

## 1. Introduction.

The optimal control of hyperbolic equations has received a lot of attention in the research literature. This is due to the large variety of applications of hyperbolic equations to physical and engineering problems.
One important direction of research in the optimal control of hyperbolic equations is the study of necessary conditions for optimality using characteristic coordinates for a controlled system governed by quasilinear second order hyperbolic equations. This has led to a family of optimal control problems for hyperbolic equations in characteristic coordinates, also known as Goursat-Darboux systems.
From a mathematical point of view, the place of optimal control for Goursat-Volterra systems in the overall framework of calculus of variations and optimal control theory can be illustrated in the following table:

| **One independent variable** | **Several independent variables** |
|---|---|
| **(1). Calculus of variations for single integrals.** | **(4). Calculus of variations for multiple integrals.** |
| **(2). Optimal control of ordinary differential equations.** | **(5). Optimal control of hyperbolic systems, in particular Roesser-type systems and Goursat-Darboux systems.** |
| **(3). Optimal control of Volterra integral equations.** | *__(6). Optimal control of Goursat-Volterra systems__ (for the first time in the present paper)**.** |

**Table 1.1.**
**Calculus of variations and optimal control in one and in several independent variables.**

The topic of optimal control of Goursat-Volterra systems is being treated for the first time in the present paper. The other topics shown in the table above have been previously treated by several authors. It is well known that in all those other topics, # (1) through (5) in the table above, there are new non-trivial developments as we move either horizontally or vertically in this table. The extremum principle of Pontryagin is generally acknowledged as the separation step of optimal control from the calculus of variations; the optimal control of Volterra integral equations has the additional feature that the Hamiltonian is a functional of the co-state, and not merely a function; the theory of calculus of variations for multiple integrals, beginning in the works of Gauss and further developed in the works of Hilbert, Carathéodory, Haar, Lichtenstein, Prange, Boerner, Kötzler, et al., has many new features that do not exist in calculus of variations for single integrals, such as Haar's lemma, Carathéodory's development of a counterpart of the Legendre transformation for multiple integral problems, etc. The transition to optimal control for hyperbolic systems presents, among other things, the new feature that, without strong differentiability assumptions on an optimal control, the Hamiltonian equations cannot be written in differential form, but only in the form of integral equations. The results in the present paper, for Goursat-Volterra systems, have the further features that (i) in addition to the case of Goursat-Darboux problems, the Hamiltonian is a functional (not just a function) of the co-state, and (ii) an extremum principle of Pontryagin's type is generally not possible, but a partial extremum principle is valid, and an extremum principle holds only in certain particular cases. These features will be made precise later in this paper.

The optimal control of Goursat-Darboux partial differential equations has been treated in, among other publications, [B1-B4, BDMO, VST]; the optimal control of ordinary Volterra equations (i.e. one independent variable) has been treated in, among other places, [B5, S, V]. The theory of optimal control of Goursat-Volterra equations is not simply a combination of Goursat-Darboux problems and ordinary Volterra equations, but rather it requires a separate treatment *ab initio*; the optimality conditions for Goursat-Volterra control problems cannot be simply guessed on the basis of Goursat-Darboux and ordinary Volterra control problems.

## 2. The basic Goursat-Volterra optimal control problem.

In this section, we shall define a certain type of hereditary systems of Goursat-Darboux type and the associated optimal control problems, and we shall examine the simplest properties of such problems.
We consider a system with partitioned state $y(s,t)=\text{col}[y_i(s,t):1\le i\le N]$ where "col" signifies a column vector obtained by placing the column vectors $y_1(s,t), y_2(s,t),...$ one underneath the other in a partioned column arrangement, and each $y_i(s,t)$ is a column vector in $\mathbb{R}^{n_i}$. We set

$$y_{(i)}(s,t):=\text{col}[y_j(s,t): 1\le j\le i-1], \text{ for } 2\le i\le N;$$
$$y_{[i]}(s,t):=\text{col}[y_j(s,t): 1\le j\le i], \text{ for } 1\le i\le N \qquad \text{--- (2.1)}$$

The controlled system will have 3 control functions $u_1(s)$, $u_2(t)$, $u_{12}(s,t)$. The system has cascaded form with respect to the components (vector components, i.e. parts of a partitioned column vector) of the state:

$$\begin{aligned}
y_i(s,t) &= f_{0,i}(s,t,y_{(i)}(s,t),u_1(s),u_2(t),u_{12}(s,t))+\\
&+\int_0^s f_{1,i}(s,t,\sigma,y_{[i]}(\sigma,t),u_1(\sigma),u_2(t),u_{12}(\sigma,t))d\sigma+\\
&+\int_0^t f_{2,i}(s,t,\tau,y_{[i]}(s,\tau),u_1(s),u_2(\tau),u_{12}(s,\tau))d\tau+\\
&+\int_0^s\int_0^t f_{12,i}(s,t,\sigma,\tau,y_{[i]}(\sigma,\tau),u_1(\sigma),u_2(\tau),u_{12}(\sigma,\tau))d\tau d\sigma
\end{aligned} \qquad \text{--- (2.2)}$$

The cost functional associated with (2.2) is

$$\begin{aligned}
J &:= F_0(y(A,B),u_1(A),u_2(B),u_{12}(A,B))+\\
&+\int_0^A F_1(s,y(s,B),u_1(s),u_2(B),u_{12}(s,B))ds+\\
&+\int_0^B F_2(t,y(A,t),u_1(A),u_2(t),u_{12}(A,t))dt+\\
&+\int_0^A\int_0^B F_{12}(s,t,y(s,t),u_1(s),u_2(t),u_{12}(s,t))dtds
\end{aligned} \qquad \text{--- (2.3)}$$

The model (2.2-2.3) encompasses the general form of controlled (non-hereditary) Goursat-Darboux systems.

The general problem of optimal control for systems governed by partial differential equations of Goursat-Darboux type has the following form of state dynamics:

$$\frac{\partial^2 x(s,t)}{\partial s \partial t} = \varphi(s,t,x(s,t),x_s(s,t),x_t(s,t),v(s,t));$$
$$x(s,0) = x_1(s);\ x(0,t) = x_2(t);$$
$$\frac{dx_1(s)}{ds} = \varphi_1(s,x_1(s),v_1(s));\ \frac{dx_2(t)}{dt} = \varphi_2(t,x_2(t),v_2(t));\ x_1(0) = x_2(0) = x_0$$

--- (2.4)

The cost functional to be minimized for (2.3) is

$$J_1 := \int_0^A \int_0^B \Phi(s,t,x(s,t),x_s(s,t),x_t(s,t),v(s,t),v_1(s),v_2(s))dtds +$$
$$+ \int_0^A \Phi_1(s,x(s,B),x_s(s,B),v_1(s))ds + \int_0^B \Phi_2(t,x(A,t),x_t(A,t),v_2(t))dt +$$
$$+ \Phi_0(x(A,B))$$

--- (2.5)

We write the equations of the group (2.4) in integral form:

$$x(s,t) = x_1(s) + x_2(t) - x_0 +$$
$$+ \int_0^s \int_0^t \varphi(\sigma,\tau,x(\sigma,\tau),x_s(\sigma,\tau),x_t(\sigma,\tau),v(\sigma,\tau))d\tau d\sigma;$$
$$x_1(s) = x_0 + \int_0^s \varphi_1(\sigma,x_1(\sigma),v_1(\sigma))d\sigma;\ x_2(t) = x_0 + \int_0^t \varphi_2(\tau,x_2(\tau),v_2(\tau))d\tau$$

--- (2.6)

We set

$$p(s,t) := x_s(s,t),\ q(s,t) := x_t(s,t)$$

--- (2.7)

Then the first equation in (2.6) becomes

$$x(s,t) = x_1(s) + x_2(t) - x_0 + \int_0^s \int_0^t \varphi(\sigma,\tau,x(\sigma,\tau),p(\sigma,\tau),q(\sigma,\tau),v(\sigma,\tau))d\tau d\sigma$$

--- (2.8)

Differentiation of (2.8) with respect to s and t yields

$$p(s,t) = \varphi_1(s, x_1(s), v_1(s)) + \int_0^t \varphi(s,\tau, x(s,\tau), p(s,\tau), q(s,\tau), v(s,\tau))d\tau;$$

$$q(s,t) = \varphi_2(t, x_2(t), v_2(t)) + \int_0^s \varphi(\sigma, t, x(\sigma,t), p(\sigma,t), q(\sigma,t), v(\sigma,t))d\sigma$$

--- (2.9)

We define

$$y_1(s,t) := \begin{bmatrix} x_1(s) \\ x_2(t) \end{bmatrix}, \ y_2(s,t) := \begin{bmatrix} x(s,t) \\ p(s,t) \\ q(s,t) \end{bmatrix}$$

--- (2.10)

Then the system of (2.6) through (2.9) is seen to be a particular case of (2.2). Also, the functional $J_1$ becomes a particular case of J.

Many of the features of the system consisting of (2.2) and (2.3) are contained in the simpler problem consisting of the state dynamics equations in the simplified form (2.11) below and the problem of minimizing the same functional J as in (2.3) but with the new simpler equations of state:

$$y(s,t) = f_0(s,t,u_1(s),u_2(t),u_{12}(s,t)) +$$
$$+\int_0^s f_1(s,t,\sigma, y(\sigma,t), u_1(\sigma), u_2(t), u_{12}(\sigma,t))d\sigma +$$
$$+\int_0^t f_2(s,t,\tau, y(s,\tau), u_1(s), u_2(\tau), u_{12}(s,\tau))d\tau +$$
$$+\int_0^s \int_0^t f_{12}(s,t,\sigma,\tau, y(\sigma,\tau), u_1(\sigma), u_2(\tau), u_{12}(\sigma,\tau))d\tau d\sigma$$

--- (2.11)

where now the state y(.) takes values in $\mathbb{R}^n$.

In some cases, we shall give complete proofs only for the model (2.11).

Next, we outline the existence and uniqueness theory for integral systems of the form (2.11).

Exclusively in the context of discussing existence and uniqueness questions, we make the following assumptions concerning the control functions and the functions $f_0, f_1, f_2, f_{12}$:

(i). $u_1(.), u_2(.), u_{12}(.,.)$ are bounded measurable functions taking values in the compact sets $U_1, U_2, U_{12}$, respectively.

(ii). The functions $f_0, f_1, f_2, f_{12}$ are bounded measurable functions with respect to $s, t, \sigma, \tau$ for all $(s, \sigma, t, \tau)$ such that $0 \le \sigma \le s \le a, 0 \le \tau \le t \le b$, and continuous with respect to $y, u_1, u_2, u_{12}$.

(iii). The functions $f_1, f_2, f_{12}$ are Lipschitz with respect to y, with Lipshitz constants $L_1, L_2, L_{12}$, respectively. ///

Under these conditions, the system (2.11) has a unique solution in $L^\infty([0,A]\times[0,B])$. This is proved by a standard contraction mapping argument, similar to the case of single Volterra equations. For completeness, we briefly outline one proving that. We define an operator S on $L^\infty([0,A]\times[0,B])$ by

$$
\begin{aligned}
(Sz)(s,t) &= f_0(s,t,u_1(s),u_2(t),u_{12}(s,t)) + \\
&+ \int_0^s f_1(s,t,\sigma,z(\sigma,t),u_1(\sigma),u_2(t),u_{12}(\sigma,t))d\sigma + \\
&+ \int_0^t f_2(s,t,\tau,z(s,\tau),u_1(s),u_2(\tau),u_{12}(s,\tau))d\tau + \\
&+ \int_0^s \int_0^t f_{12}(s,t,\sigma,\tau,z(\sigma,\tau),u_1(\sigma),u_2(\tau),u_{12}(\sigma,\tau))d\tau d\sigma.
\end{aligned}
$$

Then the problem of solving (2.11) amounts to finding a fixed point of S,

Using the weighted norm (which is equivalent to the standard norm on $L^\infty([0,a]\times[0,b])$) $\| z \|_\mu := \operatorname{ess\,sup}\{\exp(-\mu(s+t)) \,|\, z(s,t) \,| : (s,t) \in [0,a]\times[0,b]\}$, with $\mu > 0$, we have, for every two functions $z_1, z_2$

$$\exp(-\mu(s+t))\,|(Sz_1)(s,t)-(Sz_2)(s,t)|\le$$

$$\le \exp(-\mu(s+t))\left[L_1\int_0^s |z_1(\sigma,t)-z_2(\sigma,t)|d\sigma + L_2\int_0^t |z_1(s,\tau)-z_2(s,\tau)|d\tau + \right.$$

$$\left. + L_{12}\int_0^s\int_0^t |z_1(\sigma,\tau)-z_2(\sigma,\tau)|d\tau d\sigma \right]\le$$

$$\le \exp(-\mu(s+t))\left[L_1\int_0^s \exp(\mu\sigma)\|z_1-z_2\|_\mu d\sigma + L_2\int_0^t \exp(\mu\tau)\|z_1-z_2\|_\mu d\tau + \right.$$

$$\left. + L_{12}\int_0^s\int_0^t \exp(\mu(\sigma+\tau))\|z_1-z_2\|_\mu d\tau d\sigma \right]=$$

$$=\{\mu^{-1}[\exp(-\mu t)-\exp(-\mu(s+t))]L_1+\mu^{-1}[\exp(-\mu s)-\exp(-\mu(s+t))]L_2+$$

$$+\mu^{-2}[1-\exp(-\mu t)-\exp(-\mu s)+\exp(-\mu(s+t))]L_{12}\}\|z_1-z_2\|_\mu\le$$

$$\le\{\mu^{-1}[1-\exp(-\mu(a+b))]L_1+\mu^{-1}[1-\exp(-\mu(a+b))]+$$

$$+\mu^{-2}[2-\exp(-\mu a)-\exp(-\mu b)]\}\|z_1-z_2\|_\mu\,.$$

The expression

$\{\mu^{-1}[1-\exp(-\mu(A+B))]L_1+\mu^{-1}[1-\exp(-\mu(A+B))]+\mu^{-2}[2-\exp(-\mu A)-\exp(-\mu B)]\}$
goes to 0 as $\mu\to\infty$, thus it can be made $\le q$ for $\mu$ sufficiently large, if q is any predetermined number in [0, 1), and this makes the operator S a contraction on $L^\infty([0,A]\times[0,B])$ with respect to the norm $\|\cdot\|_\mu$, and then existence and uniqueness of a fixed point follows from the standard contraction mapping theorem.

We note that related existence results, for less general problems and without considering controlled systems, have been obtained by a number of authors, for example [SU].

Although the above conditions yield the existence and uniqueness of a bounded measurable solution y, they are not enough to make all terms in the cost functional J meaningful. If the terms $u_1(A),u_2(B),u_{12}(A,t),u_{12}(s,B),u_{12}(A,B)$ are absent from J, then the stated conditions suffice to make J well defined. If these terms are included in J, then we need, additionally, the following set of conditions:

(iv). The limits $\lim_{s\to A^-} u_{12}(s,\cdot)\equiv u_{12}(A,\cdot)$ and $\lim_{t\to B^-} u_{12}(\cdot,t)\equiv u_{12}(\cdot,B)$ exist in $L^\infty([0,B])$, $L^\infty([0,A])$, respectively; the limits $\lim_{s\to A^-} u_1(s)\equiv u_1(A)$ and $\lim_{t\to B^-} u_2(t)\equiv u_2(B)$ exist; the following equality involving limits holds:

$$\lim_{t\to B^-} u_{12}(A,t) = \lim_{s\to A^-} u_{12}(s,B) = \lim_{(s,t)\to(A,B)} u_{12}(s,t) \equiv u_{12}(A,B)\,.\ ///$$

An alternative, and more restrictive, set of conditions is provided below. It is possible, judging by analogy with the case of controlled ordinary differential equations, that, in cases in which optimal controls can actually be constructed, they may turn out to be piecewise continuous functions. It is, therefore, of practical interest to have optimality conditions formulated for piecewise continuous controls. The conditions below make exact the concept of "piecewise continuous controls" for the problem under consideration.

(iv*). There exist finitely many piecewise $C^1$ curves $(\gamma_j)$ in $G := [0,A]\times[0,B]$, such that the subdomains $G_i$ determined by the requirements $G_i \subseteq \operatorname{int} G\ \forall i$, each $G_i$ does not contain any part of any $(\gamma_j)$, and each $\partial G_i$ is a subset of $\bigcup_j (\gamma_j)$, are open sets. The function $u_{12}(\cdot,\cdot)$ is continuous on each $G_i$, and the restriction of $u_{12}(\cdot,\cdot)$ to $G_i$ can be extended to a function that is continuous in $\overline{G_i}$. The limits $\lim_{s\to A^-} u_{12}(s,t) \equiv u_{12}(A,t)$ and $\lim_{t\to B^-} u_{12}(s,t) \equiv u_{12}(s,B)$ exist for all except possibly finitely many values of t, respectively finitely many values of s, and they exist uniformly in t, respectively uniformly in s; the functions $u_{12}(A,\cdot)$ and $u_{12}(\cdot,B)$ are piecewise continuous (in the usual sense for functions of one variable, i.e. they can have at most a finite number of jump discontinuities with finite limits from the left and from the right at each point of discontinuity) on [0, B], respectively [0, A]; the following holds:

$$\lim_{(s,t)\to(A^-,B^-)} u_{12}(s,t) = \lim_{t\to B^-} u_{12}(A,t) = \lim_{s\to A^-} u_{12}(s,B) \equiv u_{12}(A,B)\ ;$$

$u_1(\cdot)$, $u_2(\cdot)$ are piecewise continuous on [0, A], [0, B], respectively, and the limits $\lim_{s\to A^-} u_1(s) \equiv u_1(A)$, $\lim_{t\to B^-} u_2(t) \equiv u_2(B)$ exist. ///

We conclude this section with a few remarks on applied scenarios that can lead to optimal control problems for Goursat-Volterra systems of the type we have defined above.

One category of classical applications of Goursat-Darboux and Roesser problems concerns phenomena with exchange of a quantity (e.g. mass or energy) between a moving medium and a stationary medium. Such phenomena include, among other things, gas chromatography, drying processes, air sparging methods of remediation of contaminated groundwater sites, and "moving bed operations" in chemical industry. In some of these problems, certain memoryless constitutive laws can be replaced with laws with memory, and the corresponding equations become Goursat-Volterra equations. For example, consider the problem of gas chromatography with one-dimensional flow of gas in a porous medium. If $v(s,t)$ is the speed of the flow at position s at time t, $\varphi(s,t)$ is the mass density of flowing gas, $a(s,t)$ is the mass density of absorbed gas, and $y(s,t)$ is the mass density of gas in chemical equlibrium with the absorbed gas, then the mass balance equations can be written as

$$v(s,t)\frac{\partial\varphi(s,t)}{\partial s}+\frac{\partial}{\partial t}(a(s,t)+\varphi(s,t))=0;\ \frac{\partial a(s,t)}{\partial t}=\beta(\varphi(s,t)-y(s,t)) \qquad \text{--- (2.12)}$$

where β>0 is an absorption coefficient. The system of two mass balance equations with 3 unknowns φ, a, and y, is supplemented by a constitutive law, which here we take as a law with memory,

$$a(s,t)=\int_0^t K(s,t,\tau)y(s,\tau)\,d\tau$$

(In the case of a constitutive law without memory, the corresponding equation is the so-called adsorption isotherm.)

It follows from the hereditary constitutive law that

$$a_t(s,t)=K(s,t,t)y(s,t)+\int_0^t K_t(s,t,\tau)y(s,\tau)\,d\tau$$

Then we can eliminate $a_t(s,t)$ between the above integral equation and the second equation out of the two mass balance equations:

$$[\beta+K(s,t,t)]y(s,t)=\beta\varphi(s,t)+\int_0^t K_t(s,t,\tau)y(s,\tau)\,d\tau$$

Assuming $K(s,t,t)\geq 0$, the above equation can be written in the form

$$y(s,t)=\ell_0(s,t)\varphi(s,t)+\int_0^t \frac{K_t(s,t,\tau)}{\beta+K(s,t,t)}y(s,\tau)\,d\tau;\ \ell_0(s,t):=\frac{\beta}{\beta+K(s,t,t)}$$

and, by using the resolvent kernel for a linear Volterra equation, we can represent the solution y(s, t) as

$$y(s,t)=\ell_0(s,t)\varphi(s,t)+\int_0^t \ell_1(s,t,\tau)\varphi(s,\tau)\,d\tau$$

from which it follows that

$$y_t(s,t) = \ell_{0,t}(s,t)\varphi(s,t) + \ell_0(s,t)\varphi_t(s,t) + \ell_1(s,t,t)\varphi(s,t) + \int_0^t \ell_{1,t}(s,t,\tau)\varphi(s,\tau)\,d\tau \qquad \text{--- (2.13)}$$

We differentiate eq. (2.12) with respect to t to obtain

$$v_t(s,t)\varphi_s(s,t) + v(s,t)\varphi_{st}(s,t) + \varphi_{tt}(s,t) + \beta(\varphi_t(s,t) - y_t(s,t)) = 0$$

Under the condition that the second derivative of φ with respect to time is negligible compared to the other terms, and using the expression (2.13) for $y_t(s,t)$, we can write the following equation:

$$\frac{\partial^2\varphi(s,t)}{\partial s\partial t} = a_0(s,t)\varphi(s,t) + a_1(s,t)\varphi_s(s,t) + a_2(s,t)\varphi_t(s,t) + $$
$$+ a_3(s,t)\int_0^t \ell_{1,t}(s,t,\tau)\,\varphi(s,\tau)\,d\tau \qquad \text{--- (2.14)}$$

We write (2.14) in integral form:

$$\varphi(s,t) = \varphi(s,0) + \varphi(0,t) - \varphi(0,0) + \int_0^s \int_0^t \{a_0(\sigma,\tau)\varphi(\sigma,\tau) + a_1(\sigma,\tau)\varphi_\sigma(\sigma,\tau) +$$
$$+ a_2(\sigma,\tau)\varphi_\tau(\sigma,\tau)\}\,d\tau\,d\sigma + \int_0^s \int_0^t \int_0^{\tau_1} a_3(\sigma,\tau_1)\ell_{1,\tau_1}(\sigma,\tau_1,\tau)\varphi(\sigma,\tau)\,d\tau\,d\tau_1\,d\sigma \qquad \text{--- (2.15)}$$

We define $\tilde{a}_3(t,\sigma,\tau) := \int_\tau^t a_3(\sigma,\tau_1)\ell_{1,\tau_1}(\sigma,\tau_1,\tau)\,d\tau_1$, and then eq. (2.15) can be written as

$$\varphi(s,t) = \varphi(s,0) + \varphi(0,t) - \varphi(0,0) + \int_0^s \int_0^t \{a_0(\sigma,\tau)\varphi(\sigma,\tau) + a_1(\sigma,\tau)\varphi_\sigma(\sigma,\tau) +$$
$$+ a_2(\sigma,\tau)\varphi_\tau(\sigma,\tau) + \tilde{a}_3(t,\sigma,\tau)\varphi(\sigma,\tau)\}\,d\tau\,d\sigma \qquad \text{--- (2.16)}$$

Then, by differentiating (2.16) with respect to each of s, t, we find

$$\varphi_s(s,t) = \varphi_s(s,0) + \int_0^t \{a_0(s,\tau)\varphi(s,\tau) + a_1(s,\tau)\varphi_s(s,\tau) + a_2(s,\tau)\varphi_\tau(s,\tau) +$$
$$+ \tilde{a}_3(t,s,\tau)\varphi(s,\tau)\}\, d\tau;$$
$$\varphi_t(s,t) = \varphi_t(0,t) + \int_0^s \{a_0(\sigma,t)\varphi(\sigma,t) + a_1(\sigma,t)\varphi_\sigma(\sigma,t) + a_2(\sigma,t)\varphi_t(\sigma,t) +$$
$$\tilde{a}_3(t,\sigma,t)\varphi(\sigma,t)\}\, d\tau + \int_0^s \int_0^t \tilde{a}_{3,t}(t,\sigma,\tau)\varphi(\sigma,\tau)\, d\tau\, d\sigma$$

--- (2.17)

and it is seen that the system of (2.17) is a Goursat-Volterra system, according to our definition, in the vector-valued unknown function $[\varphi(s,t) \quad \varphi_s(s,t) \quad \varphi_t(s,t)]^T$.
An optimal control problem may be formulated by taking the velocity $v(s,t)$ as the control, and asking for the minimization of a functional J given by

$$J := \int_0^A \int_0^B \{(\varphi(s,t) - \varphi_0(s,t))^2 + \mu(v(s,t))^2\}\, dt\, ds$$

i.e. we want the state $\varphi$ to be as close as possible, in the mean square sense, to a given profile $\varphi_0$, without letting the velocity become too large. This can be interpreted also as an identification problem, if the velocity field is to be estimated on the basis of observations $\varphi_0$ of the state $\varphi$.

It should be emphasized that this is just one example, among many possible applications, and it is not intended as the motivation for the type of control problems we consider in the present paper.

3. Formal derivation of the Hamiltonian equations.

The formal derivation of Hamiltonian equations relies on introducing formal penalty terms and the associated Lagrange (functional) multipliers. This formal derivation is not proof of the validity of the Hamiltonian equations, but only establishes the form that we expect the Hamiltonian equations to have. The proofs will be carried out in subsequent sections.

For the problem consisting of (2.2) and (2.3), we introduce the following penalty terms and co-state functions:

$$P_0 := \sum_i \psi_{0,i}[-y_i(A,B) + f_{0,i}(A,B,y_{(i)}(A,B),u_1(A),u_2(B),u_{12}(A,B)) +$$
$$+\int_0^A f_{1,i}(A,B,\sigma,y_{[i]}(\sigma,B),u_1(\sigma),u_2(B),u_{12}(\sigma,B))d\sigma +$$
$$+\int_0^B f_{2,i}(A,B,\tau,y_{[i]}(A,\tau),u_1(A),u_2(\tau),u_{12}(A,\tau))d\tau +$$
$$+\int_0^A \int_0^B f_{12,i}(A,B,\sigma,\tau,y(\sigma,\tau),u_1(\sigma),u_2(\tau),u_{12}(\sigma,\tau))d\tau d\sigma];$$

$$P_1 := \sum_i \int_0^A \psi_{1,i}(s)[-y_i(s,B) + f_{0,i}(s,B,y_{(i)}(s,B),u_1(s),u_2(B),u_{12}(s,B)) +$$
$$+\int_0^s f_{1,i}(s,B,\sigma,y_{[i]}(\sigma,B),u_1(\sigma),u_2(B),u_{12}(\sigma,B))d\sigma +$$
$$+\int_0^B f_{2,i}(s,B,\tau,y_{[i]}(s,\tau),u_1(s),u_2(\tau),u_{12}(s,\tau))d\tau +$$
$$+\int_0^s \int_0^B f_{12,i}(s,B,\sigma,\tau,y(\sigma,\tau),u_1(\sigma),u_2(\tau),u_{12}(\sigma,\tau))d\tau d\sigma]ds;$$

$$P_2 := \sum_i \int_0^B \psi_{2,i}(t)[-y_i(A,t) + f_{0,i}(A,t,y_{(i)}(A,t),u_1(A),u_2(t),u_{12}(A,t)) +$$
$$+\int_0^A f_{1,i}(A,t,\sigma,y_{[i]}(\sigma,t),u_1(\sigma),u_2(t),u_{12}(\sigma,t))d\sigma +$$
$$+\int_0^t f_{2,i}(A,t,\tau,y_{[i]}(A,\tau),u_1(A),u_2(\tau),u_{12}(A,\tau))d\tau +$$
$$+\int_0^A \int_0^t f_{12,i}(A,t,\sigma,\tau,y(\sigma,\tau),u_1(\sigma),u_2(\tau),u_{12}(\sigma,\tau))d\tau d\sigma]dt;$$

$$P_{12} := \sum_i \int_0^A \int_0^B \psi_{12,i}(s,t)[-y_i(s,t) + f_{0,i}(s,t,y_{(i)}(s,t),u_1(s),u_2(t),u_{12}(s,t)) +$$
$$+\int_0^s f_{1,i}(s,t,\sigma,y_{[i]}(\sigma,t),u_1(\sigma),u_2(t),u_{12}(\sigma,t))d\sigma +$$
$$+\int_0^t f_{2,i}(s,t,\tau,y_{[i]}(s,\tau),u_1(s),u_2(\tau),u_{12}(s,\tau))d\tau +$$
$$+\int_0^s \int_0^t f_{12,i}(s,t,\sigma,\tau,y_{[i]}(\sigma,\tau),u_1(\sigma),u_2(\tau),u_{12}(\sigma,\tau))d\tau d\sigma]dtds$$

--- (3.1)

The Lagrangian functional associated with (2.2) and (2.3) is

$$L := J + P_0 + P_1 + P_2 + P_{12}$$

--- (3.2)

We denote by $\delta^{\sim}L$ the variation of L obtained by using only variations $\delta^{\sim}y$ of the state y, while treating y as if it were independent from $u_1, u_2, u_{12}$.

We use the following notation: for any $\mathbb{R}^N$ -valued function g(z,w,...), where z takes values in $\mathbb{R}^M$, w takes values in $\mathbb{R}^Q$, etc., we denote by $\nabla_z g$, $\nabla_w g$,.... the Jacobians of g with respect to z, w, ..., i.e. the matices $\nabla_z g = \left[\frac{\partial g_i}{\partial z_j}\right]_{\substack{1\le i\le N\\ 1\le j\le M}}$, $\nabla_w g = \left[\frac{\partial g_i}{\partial w_j}\right]_{\substack{1\le i\le N\\ 1\le j\le Q}}$, etc. In these formulae, the subscript i denotes rows and j denotes columns of the indicated matrices. According to this notation, the gradient of a real-valued function is a row vector.

Then we have

$$\delta^{\sim}P_0 = \sum_{i} (-\psi_{0,i})\delta^{\sim}y_i(A,B) + \sum_{i,j} \psi_{0,i}\nabla_{y_j} f_{0,i}(A,B,...)\delta^{\sim}y_j(A,B) +$$

$$+\sum_{i,j} \int_0^A \psi_{0,i}\nabla_{y_j} f_{1,i}(A,B,s,...)\delta^{\sim}y_j(s,B)ds +$$

$$+\sum_{i,j} \int_0^B \psi_{0,i}\nabla_{y_j} f_{2,i}(A,B,t,...)\delta^{\sim}y_j(A,t)dt +$$

$$+\sum_{i,j} \int_0^A \int_0^B \psi_{0,i}\nabla_{y_j} f_{12,i}(A,B,s,t,...)\delta^{\sim}y_j(s,t)dtds$$

--- (3.3)

$$\delta^{\sim}P_1 = \sum_{i} \Big\{ \int_0^A \psi_{1,i}(s)[-\delta^{\sim}y_i(s,B) + \sum_{j:\, j<i} \nabla_{y_j} f_{0,i}(s,B,...)\delta^{\sim}y_j(s,B)]ds +$$

$$+\sum_{j:\, j\le i} \int_0^A \int_0^s \psi_{1,i}(s)\nabla_{y_j} f_{1,i}(s,B,\sigma,...)\delta^{\sim}y_j(\sigma,B)d\sigma ds +$$

$$+\sum_{j:\, j\le i} \int_0^A \int_0^B \psi_{1,i}(s)\nabla_{y_j} f_{2,i}(s,B,\tau,...)\delta^{\sim}y_j(s,\tau)d\tau ds +$$

$$+\sum_{j:\, j\le i} \int_0^A \int_0^s \int_0^B \psi_{1,i}(s)\nabla_{y_j} f_{12,i}(s,B,\sigma,\tau,...)\delta^{\sim}y_j(\sigma,\tau)d\tau d\sigma ds \Big\}$$

--- (3.4)

$$\delta^{\sim}P_2 = \sum_{i} \Big\{ \int_0^B \psi_{2,i}(t)[-\delta^{\sim}y_i(A,t) + \sum_{j:\, j<i} \nabla_{y_j} f_{0,i}(A,t,...)\delta^{\sim}y_j(A,t)]dt +$$

$$+\sum_{j:\, j\le i} \int_0^A \int_0^B \psi_{2,i}(t)\nabla_{y_j} f_{1,i}(A,t,\sigma,...)\delta^{\sim}y_j(\sigma,t)dtd\sigma +$$

$$+\sum_{j:\, j\le i} \int_0^B \int_0^t \psi_{2,i}(t)\nabla_{y_j} f_{2,i}(A,t,\tau,...)\delta^{\sim}y_j(A,\tau)d\tau dt +$$

$$+\sum_{j:\, j\le i} \int_0^A \int_0^B \int_0^t \psi_{2,i}(t)\nabla_{y_j} f_{12,i}(A,t,\sigma,\tau,...)\delta^{\sim}y_j(\sigma,\tau)d\tau dtd\sigma \Big\}$$

--- (3.5)

$$\delta^{\sim}P_{12} = \sum_i \left\{ \int_0^A \int_0^B \psi_{12,i}(s,t)[-\delta^{\sim}y_i(s,t) + \right.$$

$$+ \sum_{j:\, j<i} \nabla_{y_j} f_{0,i}(s,t,...)\delta^{\sim}y_j(s,t)]dtds +$$

$$+ \sum_{j:\, j\leq i} \int_0^A \int_0^B \int_0^s \psi_{12,i}(s,t)\nabla_{y_j} f_{1,i}(s,t,\sigma,...)\delta^{\sim}y_j(\sigma,t)d\sigma dtds +$$

$$+ \sum_{j:\, j\leq i} \int_0^A \int_0^B \int_0^t \psi_{12,i}(s,t)\nabla_{y_j} f_{2,i}(s,t,\tau,...)\delta^{\sim}y_j(s,\tau)d\tau dtds +$$

$$\left. + \sum_{j:\, j\leq i} \int_0^A \int_0^B \int_0^s \int_0^t \psi_{12,i}(s,t)\nabla_{y_j} f_{12,i}(s,t,\sigma,...)\delta^{\sim}y_j(\sigma,\tau)d\tau d\sigma dtds \right\}$$

--- (3.6)

$$\delta^{\sim}J = \sum_i \left\{ \psi_{0,i}\nabla_{y_i} F_0(y(A,B),...)\delta^{\sim}y_i(A,B) + \right.$$

$$+ \int_0^A \psi_{1,i}(s)\nabla_{y_i} F_1(s,...)\delta^{\sim}y_i(s,B)ds + \int_0^B \psi_{2,i}(t)\nabla_{y_i} F_2(t,...)\delta^{\sim}y_i(A,t)dt +$$

$$\left. + \int_0^A \int_0^B \psi_{12,i}(s,t)\nabla_{y_i} F_{12}(s,t,...)\delta^{\sim}y_i(s,t)dtds \right\}$$

(3.7)

By setting the variation of L equal to zero, and taking into account the above formulae for the variations of the individual terms that make up L, and after performing a few changes in the order of integration and renaming certain variables, we obtain the following integral equations for the co-state $[\psi_0 \ \psi_1(s) \ \psi_2(t) \ \psi_{12}(s,t)]$:

$$\psi_{12,i}(s,t) = \nabla_{y_i} F_{12}(s,t,...) + \sum_{j:\, j>i} \psi_{12,j}(s,t)\nabla_{y_i} f_{0,j}(s,t,...) +$$

$$+ \sum_{j:\, j\geq i} \Bigg\{ \psi_{1,j}(s)\nabla_{y_i} f_{2,j}(s,B,t,...) + \psi_{2,j}(t)\nabla_{y_i} f_{1,j}(A,t,s,...) + \psi_{0,j}\nabla_{y_i} f_{12}(A,B,s,t,...) +$$

$$+ \int_s^A [\psi_{12,j}(\sigma,t)\nabla_{y_i} f_{1,j}(\sigma,t,s,...) + \psi_{1,j}(\sigma)\nabla_{y_i} f_{12,j}(\sigma,B,s,t,...)]\, d\sigma +$$

$$+ \int_t^B [\psi_{12,j}(s,\tau)\nabla_{y_i} f_{2,j}(s,\tau,t,...) + \psi_{2,j}(\tau)\nabla_{y_i} f_{12,j}(A,\tau,s,t,...)]\, d\tau +$$

$$+ \int_s^A \int_t^B \psi_{12,j}(\sigma,\tau)\nabla_{y_i} f_{12}(\sigma,\tau,s,t,...)\, d\tau d\sigma \Bigg\}$$

--- (3.8)

$$\psi_{1,i}(s) = \nabla_{y_i} F_1(s,...) + \sum_{j:\, j>i} \psi_{1,j}(s)\nabla_{y_i} f_{0,j}(s,B,...) +$$

$$+ \sum_{j:\, j\geq i} \Bigg\{ \psi_{0,j}\nabla_{y_i} f_{1,j}(A,B,s,...) + \int_s^A \psi_{1,j}(\sigma)\nabla_{y_i} f_{1,j}(\sigma,B,s,...)\, d\sigma \Bigg\}$$

--- (3.9)

$$\psi_{2,i}(t) = \nabla_{y_i} F_2(t,...) + \sum_{j:\, j>i} \psi_{2,j}(t)\nabla_{y_i} f_{0,j}(A,t,...) +$$

$$+ \sum_{j:\, j\geq i} \Bigg\{ \psi_{0,j}\nabla_{y_i} f_{2,j}(A,B,t,...) + \int_t^B \psi_{2,j}(\tau)\nabla_{y_i} f_{2,j}(A,\tau,t,...)\, d\tau \Bigg\}$$

--- (3.10)

$$\psi_{0,i} = \nabla_{y_i} F_0(y(A,B),...) + \sum_{j:\, j>i} \psi_{0,j} f_{0,j}(A,B,...)$$

--- (3.11)

This leads us to define 4 Hamiltonians $H_0, H_1, H_2, H_{12}$. In addition to the variables $u_1, u_2, u_{12}$, these Hamiltonians also depend on the additional variables $u_{1,1}$, $u_{2,2}$, $u_{12,1}$, $u_{12,2}$, $u_{12,12}$, where $u_{1,1}$ plays the role of $u_1(A)$, $u_{2,2}$ that of $u_2(B)$, $u_{12,1}$ that of $u_{12}(A,t)$, $u_{12,2}$ that of $u_{12}(s,B)$, and $u_{12,12}$ that of $u_{12}(A,B)$. Also, these Hamiltonians depend on the variables $y, y_{\{1\}}, y_{\{2\}}, y_{\{12\}}$ , where $y_{\{1\}}$ plays the role of $y(A,t)$, $y_{\{2\}}$ that of $y(s,B)$, and $y_{\{12\}}$ that of $y(A,B)$. The symbols $y_{\{1\},j}$, $y_{\{1\},(j)}, y_{\{1\},[j]}$, etc. will have the analogous meaning with $y_j$, $y_{(j)}$, $y_{[j]}$, for example $y_{\{1\},(j)}$ plays the role of $y_{(j)}(A,t)$, and so forth.

We define

$$H_0(y_{\{12\}}, u_{1,1}, u_{2,2}, u_{12,12}, \psi_0) :=$$

$$= F_0(y_{\{12\}}, u_{1,1}, u_{2,2}, u_{12,12}) + \sum_j \psi_{0,j} f_{0,j}(A, B, y_{\{12\},(j)}, u_{1,1}, u_{2,2}, u_{12,12});$$

$$H_1(s, y_{\{2\}}, u_1, u_{2,2}, u_{12,2}, \psi_0, \psi_1(\cdot)) := F_1(s, y_{\{2\}}, u_1, u_{2,2}, u_{12,2}) +$$

$$+ \sum_j \Bigg[ \psi_{1,j}(s) f_{0,j}(s, B, y_{\{2\},(j)}, u_1, u_{2,2}, u_{12,2}) + \psi_{0,j} f_{1,j}(A, B, s, y_{\{2\},[j]}, u_1, u_{2,2}, u_{12,2}) +$$

$$+ \int_s^A \psi_{1,j}(\sigma) f_{1,j}(\sigma, B, s, y_{\{2\},[j]}, u_1, u_{2,2}, u_{12,2})\, d\sigma \Bigg];$$

$$H_2(t, y_{\{1\}}, u_{1,1}, u_2, u_{12,1}, \psi_0, \psi_2(\cdot)) := F_2(t, y_{\{1\}}, u_{1,1}, u_2, u_{12,1}) +$$

$$+ \sum_j \Bigg[ \psi_{2,j}(t) f_{0,j}(A, t, y_{\{1\},(j)}, u_{1,1}, u_2, u_{12,1}) + \psi_{0,j} f_{2,j}(A, B, t, y_{\{1\},[j]}, u_{1,1}, u_2, u_{12,1}) +$$

$$+ \int_t^B \psi_{2,j}(\tau) f_{2,j}(A, \tau, t, y_{\{1\},[j]}, u_{1,1}, u_2, u_{12,1})\, d\tau \Bigg];$$

$$H_{12}(s, t, y, u_1, u_2, u_{12}, \psi_0, \psi_1(.), \psi_2(.), \psi_{12}(.,.)) := F_{12}(s, t, y, u_1, u_2, u_{12}) +$$

$$+ \sum_j \Bigg[ \psi_{0,j} f_{12,j}(A, B, s, t, y_{[i]}, u_1, u_2, u_{12}) + \psi_{1,j}(s) f_{2,j}(s, B, t, y_{[i]}, u_1, u_2, u_{12}) +$$

$$+ \psi_{2,j}(t) f_{1,j}(A, t, s, y_{[i]}, u_1, u_2, u_{12}) + \psi_{12,j}(s,t) f_{0,j}(s, t, y_{(j)}, u_1, u_2, u_{12}) +$$

$$+ \int_s^A [\psi_{12,j}(\sigma, t) f_{1,j}(\sigma, t, s, y_{[i]}, u_1, u_2, u_{12}) + \psi_{1,j}(\sigma) f_{12,j}(\sigma, B, s, t, y_{[i]}, u_1, u_2, u_{12})]\, d\sigma +$$

$$+ \int_t^B [\psi_{12,j}(s, \tau) f_{2,j}(s, \tau, t, y_{[i]}, u_1, u_2, u_{12}) + \psi_{2,j}(\tau) f_{12,j}(A, \tau, s, t, y_{[i]}, u_1, u_2, u_{12})]\, d\tau +$$

$$+ \int_s^A \int_t^B \psi_{12,j}(\sigma, \tau) f_{12,j}(\sigma, \tau, s, t, y_{[i]}, u_1, u_2, u_{12})\, d\tau\, d\sigma \Bigg]$$

--- (3.12)

The Hamiltonian equations for the co-state can now be written in the form

$$\psi_{0,i} = \nabla_{y_{\{12\},i}} H_0(y_{\{12\}}, u_{1,1}, u_{2,2}, u_{12,12}, \psi_0);$$
$$\psi_{1,i}(s) = \nabla_{y_{\{2\},i}} H_1(s, y_{\{2\}}, u_1, u_{2,2}, u_{12,2}, \psi_0, \psi_1(\cdot));$$
$$\psi_{2,i}(t) = \nabla_{y_{\{1\},i}} H_2(t, y_{\{1\}}, u_{1,1}, u_2, u_{12,1}, \psi_0, \psi_2(\cdot));$$
$$\psi_{12,i}(s,t) = \nabla_{y_i} H_{12}(s, t, y, u_1, u_2, u_{12}, \psi_0, \psi_1(.), \psi_2(.), \psi_{12}(.,.))$$

--- (3.13)

For later reference, we also record below the Hamiltonian equations for the simplified model (2.11):

$$\psi_{12}(s,t) = \nabla_y F_{12}(s,t,...) + \psi_1(s)\nabla_y f_2(s,B,t,...) + \psi_2(t)\nabla_y f_1(A,t,s,...) +$$
$$+ \psi_0 \nabla_y f_{12}(A,B,s,t,...) + \int_s^A [\psi_{12}(\sigma,t)\nabla_y f_1(\sigma,t,s,...) +$$
$$+ \psi_1(\sigma)\nabla_y f_{12}(\sigma,B,s,t,...)]d\sigma +$$
$$+ \int_t^B [\psi_{12}(s,\tau)\nabla_y f_2(s,\tau,t,...) + \psi_2(\tau)\nabla_y f_{12}(A,\tau,s,t,...)]d\tau +$$
$$+ \int_s^A \int_t^B \psi_{12}(\sigma,\tau)\nabla_y f_{12}(\sigma,\tau,s,t,...)d\tau d\sigma;$$
$$\psi_1(s) = \nabla_y F_1(s,...) + \psi_0 \nabla_y f_1(A,B,s,...) + \int_s^A \psi_1(\sigma)\nabla_y f_1(\sigma,B,s,...)d\sigma;$$
$$\psi_2(t) = \nabla_y F_2(t,...) + \psi_0 \nabla_y f_2(A,B,t,...) + \int_t^B \psi_2(\tau)\nabla_y f_2(A,\tau,t,...)d\tau;$$
$$\psi_0 = \nabla_y F_0(y(A,B),...)$$

--- (3.14)

The Hamiltonians for the simplified model (2.11) are

$$h_0(y_{\{12\}}, u_{1,1}, u_{2,2}, u_{12,12}) := F_0(y_{\{12\}}, u_{1,1}, u_{2,2}, u_{12,12});$$

$$h_1(s, y_{\{2\}}, u_1, u_{2,2}, u_{12,2}, \psi_0, \psi_1(\cdot)) := F_1(s, y_{\{2\}}, u_1, u_{2,2}, u_{12,2}) +$$
$$+ \psi_0 f_1(A, B, s, y_{\{2\}}, u_1, u_{2,2}, u_{12,2}) +$$
$$+ \int_s^A \psi_1(\sigma) f_1(\sigma, B, s, y_{\{2\}}, u_1, u_{2,2}, u_{12,2})\, d\sigma;$$

$$h_2(t, y_{\{1\}}, u_{1,1}, u_2, u_{12,1}, \psi_0, \psi_2(\cdot)) := F_2(t, y_{\{1\}}, u_{1,1}, u_2, u_{12,1}) +$$
$$+ \psi_0 f_2(A, B, t, y_{\{1\}}, u_{1,1}, u_2, u_{12,1}) + \int_t^B \psi_2(\tau) f_2(A, \tau, t, y_{\{1\}}, u_{1,1}, u_2, u_{12,1})\, d\tau;$$

$$h_{12}(s, t, y, u_1, u_2, u_{12}, \psi_0, \psi_1(.), \psi_2(.), \psi_{12}(.,.)) :=$$
$$= F_{12}(s, t, y, u_1, u_2, u_{12}) + \psi_0 f_{12}(A, B, s, t, y, u_1, u_2, u_{12}) +$$
$$+ \psi_1(s) f_2(s, B, t, y, u_1, u_2, u_{12}) + \psi_2(t) f_1(A, t, s, y, u_1, u_2, u_{12}) +$$
$$+ \int_s^A [\psi_{12}(\sigma, t) f_1(\sigma, t, s, y, u_1, u_2, u_{12}) + \psi_1(\sigma) f_{12}(\sigma, B, s, t, y, u_1, u_2, u_{12})]\, d\sigma +$$
$$+ \int_t^B [\psi_{12}(s, \tau) f_2(s, \tau, t, y, u_1, u_2, u_{12}) + \psi_2(\tau) f_{12}(A, \tau, s, t, y, u_1, u_2, u_{12})]\, d\tau +$$
$$+ \int_s^A \int_t^B \psi_{12}(\sigma, \tau) f_{12}(\sigma, \tau, s, t, y, u_1, u_2, u_{12})\, d\tau\, d\sigma$$

--- (3.15)

4. The variational equations.

In this section we shall carry out the exact calculation of the variation of the functional J of section 2. In order to somewhat reduce the notational complexity, we shall examine in detail the simplified model
Assuming adequate continuous differentiability of all functions involved, we calculate the variation of the functional J under an admissible variation $\delta u$ of the control function u. We are using the standard definition of admissible variations, i.e. a variation $\delta u(\cdot)$ is said to be admissible relative to an admissible control function $u(\cdot)$ if there exists an $\varepsilon_0 > 0$ such that, for all $\varepsilon \in [0, \varepsilon_0]$, the function $u(\cdot) + \varepsilon \delta u(\cdot)$ is an admissible control function. The admissible control functions were defined in section 2. The set of admissible variations $\delta u(\cdot)$ relative to an admissible control function $u(\cdot)$ will be denoted by $\mathbf{A}(u(\cdot))$. The variation $\delta y(s,t) = \delta y(s,t;u(.),\delta u(.))$ of the state is defined as

$$\delta y(s,t) = \lim_{\varepsilon \to 0^+} \frac{1}{\varepsilon}\left(y(s,t;u(.)+\varepsilon\delta u(.)) - y(s,t;u(.))\right) \qquad \text{--- (4.1)}$$

Similarly, the variation $\delta J = \delta J(u(.),\delta u(.))$ of the functional J is defined as

$$\delta J = \lim_{\varepsilon \to 0^+} \frac{1}{\varepsilon}\left(J(u(.)+\varepsilon\delta u(.)) - J(u(.))\right) \qquad \text{--- (4.2)}$$

Then we have

$$\begin{aligned}
\delta y(s,t) = \int_0^s & \{(\nabla_y f_1(s,t,\sigma,y(\sigma,t),u(\sigma,t)))\delta y(\sigma,t) + \\
&+ (\nabla_u f_1(s,t,\sigma,y(\sigma,t),u(\sigma,t)))\delta u(\sigma,t)\}d\sigma + \\
&+ \int_0^t \{(\nabla_y f_2(s,t,\tau,y(s,\tau),u(s,\tau)))\delta y(s,\tau) + \\
&+ (\nabla_u f_2(s,t,\tau,y(s,\tau),u(s,\tau)))\delta u(s,\tau)\}d\tau + \\
&+ \int_0^s \int_0^t \{(\nabla_y f_{12}(s,t,\sigma,\tau,y(\sigma,\tau),u(\sigma,\tau)))\delta y(\sigma,\tau) + \\
&(\nabla_u f_{12}(s,t,\sigma,y(\sigma,\tau),u(\sigma,\tau)))\delta u(\sigma,\tau)\}d\tau d\sigma
\end{aligned} \qquad \text{--- (4.3)}$$

and

$$\delta J = (\nabla_y F_0(y(A,B)))\delta y(A,B) +$$
$$+\int_0^A \{(\nabla_y F_1(s,y(s,B),u(s,B)))\delta y(s,B) +$$
$$+(\nabla_u F_1(s,y(s,B),u(s,B)))\delta u(s,B)\}ds +$$
$$+\int_0^B \{(\nabla_y F_2(t,y(A,t),u(A,t)))\delta y(A,t) +$$
$$+(\nabla_u F_2(t,y(A,t),u(A,t)))\delta u(A,t)\}dt +$$
$$+\int_0^A \int_0^B \{(\nabla_y F_{12}(s,t,y(s,t),u(s,t))\delta y(s,t) +$$
$$+\{(\nabla_y F_{12}(s,t,y(s,t),u(s,t))\delta y(s,t)\}dtds$$

--- (4.4)

As it is well known, the necessary conditions for optimality are

$$\delta J(u^*(.),\delta u(.)) \geq 0 \ \forall \delta u(.) \in \mathbf{A}(u^*(.))$$

--- (4.5)

In order to proceed further, we need a theory of resolvent kernels and duality for linear Goursat-Volterra equations of the type of (4.3). This theory is developed in the next section.

5. Linear Goursat-Volterra equations.

We consider a linear Goursat-Volterra equation

$$z(s,t)=z_0(s,t)+\int_0^s K_1(s,t,\sigma)z(\sigma,t)d\sigma+\int_0^t K_2(s,t,\tau)z(s,\tau)d\tau+$$
$$+\int_0^s\int_0^t K_{12}(s,t,\sigma,\tau)z(\sigma,\tau)d\tau d\sigma$$

--- (5.1)

For the purpose of developing a theory of resolvent kernels, we denote by K the triple $(K_1,K_2,K_{12})$, and by **K** the integral operator

$$(\mathbf{K}z)(s,t):=\int_0^s K_1(s,t,\sigma)z(\sigma,t)d\sigma+\int_0^t K_2(s,t,\tau)z(s,\tau)d\tau+$$
$$+\int_0^s\int_0^t K_{12}(s,t,\sigma,\tau)z(\sigma,\tau)d\tau d\sigma$$

--- (5.2)

We also consider a second operator **L** of the same type as **K** , and the triple $L=(L_1,L_2,L_{12})$ associated with the operator **L**:

$$(\mathbf{L}w)(s,t):=\int_0^s L_1(s,t,\sigma)w(\sigma,t)d\sigma+\int_0^t L_2(s,t,\tau)w(s,\tau)d\tau+$$
$$+\int_0^s\int_0^t L_{12}(s,t,\sigma,\tau)w(\sigma,\tau)d\tau d\sigma$$

--- (5.3)

We are interested in the composition $\mathbf{M}:=\mathbf{L}\circ\mathbf{K}$. As it will be shown, this composition is again a linear Goursat-Volterra operator, and it has the corresponding triple of kernels $M=(M_1,M_2,M_{12})$.
We have

$$
\begin{aligned}
((\mathbf{L}\circ\mathbf{K})z)(s,t) &= \int_0^s\int_0^{\sigma_1} L_1(s,t,\sigma_1)K_1(\sigma_1,t,\sigma)z(\sigma,t)d\sigma d\sigma_1 + \\
&+\int_0^s\int_0^{\sigma_1} L_1(s,t,\sigma_1)K_2(\sigma_1,t,\tau)z(\sigma_1,\tau)d\tau d\sigma_1 + \\
&+\int_0^s\int_0^{\sigma_1}\int_0^t L_1(s,t,\sigma_1)K_{12}(\sigma_1,t,\sigma,\tau)z(\sigma,\tau)d\tau d\sigma d\sigma_1 + \\
&+\int_0^t\int_0^s L_2(s,t,\tau_1)K_1(s,\tau_1,\sigma)z(\sigma,\tau_1)d\sigma d\tau_1 + \\
&+\int_0^t\int_0^{\tau_1} L_2(s,t,\tau_1)K_2(s,\tau_1,\tau)z(s,\tau)d\tau d\tau_1 + \\
&+\int_0^t\int_0^s\int_0^{\tau_1} L_2(s,t,\tau_1)K_{12}(s,\tau_1,\sigma,\tau)z(\sigma,\tau)d\tau d\sigma d\tau_1 + \\
&+\int_0^s\int_0^t\int_0^{\sigma_1} L_{12}(s,t,\sigma_1,\tau_1)K_1(\sigma_1,\tau_1,\sigma)z(\sigma,\tau_1)d\sigma d\tau_1 d\sigma_1 + \\
&+\int_0^s\int_0^t\int_0^{\tau_1} L_{12}(s,t,\sigma_1,\tau_1)K_2(\sigma_1,\tau_1,\tau)z(\sigma_1,\tau)d\tau d\tau_1 d\sigma_1 + \\
&+\int_0^s\int_0^t\int_0^{\sigma_1}\int_0^{\tau_1} L_{12}(s,t,\sigma_1,\tau_1)K_{12}(\sigma_1,\tau_1,\sigma,\tau)z(\sigma,\tau)d\tau d\sigma d\tau_1 d\sigma_1
\end{aligned}
$$

--- (5.4)

We define the 9 convolutions $\otimes_{\alpha,\beta}$, where each of $\alpha,\beta$ takes values in {1, 2, 12}, by

$$
\begin{aligned}
&(L_1 \otimes_{1,1} K_1)(s,t,\sigma) := \int_\sigma^s L_1(s,t,\sigma_1)K_1(\sigma_1,t,\sigma)d\sigma_1 ; \\
&(L_2 \otimes_{2,2} K_2)(s,t,\tau) := \int_\tau^t L_2(s,t,\tau_1)K_2(s,\tau_1,\tau)d\tau_1 ; \\
&(L_1 \otimes_{1,2} K_2)(s,t,\sigma,\tau) := L_1(s,t,\sigma)K_2(s,t,\tau) ; \\
&(L_2 \otimes_{2,1} K_1)(s,t,\sigma,\tau) := L_2(s,t,\tau)K_1(s,t,\sigma) ; \\
&(L_1 \otimes_{1,12} K_{12})(s,t,\sigma,\tau) := \int_\sigma^s L_1(s,t,\sigma_1)K_{12}(\sigma_1,t,\sigma,\tau)d\sigma_1 ; \\
&(L_{12} \otimes_{12,1} K_1)(s,t,\sigma,\tau) := \int_\sigma^s L_{12}(s,t,\sigma_1,\tau)K_1(\sigma_1,\tau,\sigma)d\sigma_1 ; \\
&(L_2 \otimes_{2,12} K_{12})(s,t,\sigma,\tau) := \int_\tau^t L_2(s,t,\tau_1)K_{12}(s,\tau_1,\sigma,\tau)d\tau_1 ; \\
&(L_{12} \otimes_{12,2} K_2)(s,t,\sigma,\tau) := \int_\tau^t L_{12}(s,t,\sigma,\tau_1)K_2(\sigma,\tau_1,\tau)d\tau_1 ; \\
&(L_{12} \otimes_{12,12} K_{12})(s,t,\sigma,\tau) := \int_\sigma^s\int_\tau^t L_{12}(s,t,\sigma_1,\tau_1)K_{12}(\sigma_1,\tau_1,\sigma,\tau)d\tau_1 d\sigma_1
\end{aligned}
$$

--- (5.5)

It follows from (5.4) and (5.5) that

$$M_1 = L_1 \otimes_{1,1} K_1,\ M_2 = L_2 \otimes_{2,2} K_2,\ M_{12} = \sum_{\alpha,\beta:\ (1,1)\neq(\alpha,\beta)\neq(2,2)} L_\alpha \otimes_{\alpha,\beta} K_\beta \quad \text{--- (5.6)}$$

We shall write

$$M = L \otimes K \quad \text{--- (5.7)}$$

This convolution is associative but, in general, noncommutative. We set

$$K^{\otimes N} := \underbrace{K \otimes K \otimes \cdots \otimes K}_{N \text{ times}} \quad \text{--- (5.8)}$$

For consistency of notation, we define the convolutions

$$\begin{aligned}
(K_1 \otimes_{1,0} z)(s,t) &:= \int_0^s K_1(s,t,\sigma) z(\sigma,t) d\sigma;\\
(K_2 \otimes_{2,0} z)(s,t) &:= \int_0^t K_2(s,t,\tau) z(s,\tau) d\tau;\\
(K_{12} \otimes_{12,0} z)(s,t) &:= \int_0^s \int_0^t K_{12}(s,t,\sigma,\tau) z(\sigma,\tau) d\tau d\sigma;\\
(K \otimes_0 z)(s,t) &:= (K_1 \otimes_{1,0} z)(s,t) + (K_2 \otimes_{2,0} z)(s,t) + (K_{12} \otimes_{12,0} z)(s,t)
\end{aligned} \quad \text{--- (5.9)}$$

The construction of the convolution $\otimes$ implies

$$L \otimes_0 (K \otimes_0 z) = (L \otimes K) \otimes_0 z \quad \text{--- (5.10)}$$

The solution of the integral equation (5.1) by the method of Picard iterations (as in section 2) leads therefore to the following formula for the N-th iteration:

$$z_{(N)}(s,t) = z_0(s,t) + \sum_{k=1}^{N} (K^{\otimes k} \otimes_0 z)(s,t) \quad \text{--- (5.11)}$$

We set

$$R(s,t,\sigma,\tau) = (R_1(s,t,\sigma), R_2(s,t,\tau), R_{12}(s,t,\sigma,\tau)) :=$$
$$= \left( \left(\sum_{k=1}^{\infty} K^{\otimes k}\right)_1 (s,t,\sigma), \left(\sum_{k=1}^{\infty} K^{\otimes k}\right)_2 (s,t,\tau), \left(\sum_{k=1}^{\infty} K^{\otimes k}\right)_{12} (s,t,\sigma,\tau) \right)$$

--- (5.12)

Under the condition that all kernels $K_1, K_2, K_{12}$ are measurable and bounded for all relevant values of $s,t,\sigma,\tau$, we can prove the convergence of the 3 series in (5.12). We set

$$\mathbf{D} := \{(s,t,\sigma,\tau) : 0 \le \sigma \le s \le A, 0 \le \tau \le t \le B\}$$

--- (5.13)

and we denote by C a common bound on $K_1, K_2, K_{12}$:

$$|K_1(s,t,\sigma)| \le C, \; |K_2(s,t,\tau)| \le C, \; |K_{12}(s,t,\sigma,\tau)| \le C, \; \forall (s,t,\sigma,\tau) \in \mathbf{D}$$

--- (5.14)

We have:

<u>Lemma 5.1.</u> The terms $(K^{\otimes k})_1$, $(K^{\otimes k})_2$ satisfy

$$|(K^{\otimes k})_1(s,t,\sigma)| \le C^k \frac{(s-\sigma)^{k-1}}{(k-1)!}, \; |(K^{\otimes k})_2(s,t,\tau)| \le C^k \frac{(t-\tau)^{k-1}}{(k-1)!}, \; \forall (s,t,\sigma,\tau) \in \mathbf{D}$$

--- (5.15)

<u>Proof:</u> According to (5.5) and (5.6), we have

$$(K^{\otimes k})_1 = (K_1)^{\otimes_1 k} \equiv \underbrace{K_1 \otimes_1 K_1 \otimes_1 \cdots \otimes_1 K_1}_{k \text{ times}}$$

--- (5.16)

and consequently the estimate $|(K^{\otimes k})_1(s,t,\sigma)| \le C^k \frac{(s-\sigma)^{k-1}}{(k-1)!}$ follows by the same method as for standard Volterra equations (equations with single integrals). The estimate for $(K^{\otimes k})_2$ is, of course, justified in the same way. ///

<u>Lemma 5.2.</u> We have the estimate

$$|(K^{\otimes k})_{12}(s,t,\sigma,\tau)| \le C^k \sum_{\lambda=1}^{N_k} \varphi_{k,\lambda}(s-\sigma,t-\tau);\ N_k := 3^k - 2$$

--- (5.17)

where each $\varphi_{k,\lambda}$ is a term of the form

$$\varphi_{k,\lambda}(s-\sigma,t-\tau) = \frac{(s-\sigma)^{\alpha_{k,\lambda}}(t-\tau)^{\beta_{k,\lambda}}}{(\alpha_{k,\lambda})!(\beta_{k,\lambda})!}$$

--- (5.18)

with $\alpha_{k,\lambda}$, $\beta_{k,\lambda}$ nonnegative integers satisfying

$$\alpha_{k,\lambda} + \beta_{k,\lambda} \ge k-2,\ \max(\alpha_{k,\lambda},\beta_{k,\lambda}) \le k-1$$

--- (5.19)

(The terms $\varphi_{k,\lambda}$ appearing in the summation in (5.17) are not necessarily distinct: the summation may include repeated terms, but each term is of the form specified in (5.18) and (5.19). The set of pairs ( $\alpha_{k,\lambda}$ , $\beta_{k,\lambda}$ ) that actually appear in the terms in the summation in (5.17) does not necessarily comprise all possible pairs of nonnegative integers that satisfy (5.19).)

Proof: The numbers $N_k$ are solutions of the finite-difference equation $N_{k+1} = 4 + 3N_k$ with initial condition $N_1 = 1$. The justification for these numbers comes from the number of terms (each term being of the form $C^k \varphi_{k,\lambda}(s-\sigma,t-\tau)$) involved in the estimation of $(K^{\otimes(k+1)})_{12}$ in each new iteration, assuming $K^{\otimes k}$ satisfies the asserted estimates. The kernel $(K^{\otimes(k+1)})_{12}$ is expressed as

$$\begin{aligned}(K^{\otimes(k+1)})_{12} &= (K_1{}^{\otimes_1 k})K_2 + (K_2{}^{\otimes_2 k})K_1 + (K_1{}^{\otimes_1 k}) \otimes_{1,12} K_{12} + \\ &+ (K_2{}^{\otimes_2 k}) \otimes_{2,12} K_{12} + (K^{\otimes k})_{12} \otimes_{12,1} K_1 + (K^{\otimes k})_{12} \otimes_{12,2} K_2 + \\ &+ (K^{\otimes k})_{12} \otimes_{12,12} K_{12}\end{aligned}$$

--- (5.20)

Now, the right-hand side of (5.20) contains the four terms
$(K_1{}^{\otimes_1 k})K_2$, $(K_2{}^{\otimes_2 k})K_1$, $(K_1{}^{\otimes_1 k}) \otimes_{1,12} K_{12}$, $(K_2{}^{\otimes_2 k}) \otimes_{2,12} K_{12}$,

which contribute, respectively, terms of the form

$$C^{k+1}\frac{(s-\sigma)^{k-1}}{(k-1)!},\ C^{k+1}\frac{(t-\tau)^{k-1}}{(k-1)!},\ C^{k+1}\frac{(s-\sigma)^{k}}{k!},\ C^{k+1}\frac{(t-\tau)^{k}}{k!}$$

to the estimation of $|(K^{\otimes(k+1)})_{12}(s,t,\sigma,\tau)|$; each of the remaining terms $(K^{\otimes k})_{12}\otimes_{12,1}K_1$, $(K^{\otimes k})_{12}\otimes_{12,2}K_2$, $(K^{\otimes k})_{12}\otimes_{12,12}K_{12}$ contributes $N_k$ terms to the estiamtion of $|(K^{\otimes(k+1)})_{12}(s,t,\sigma,\tau)|$; thus $N_{k+1}=4+3N_k$. For each term of the form $C^k\frac{(s-\sigma)^{\alpha_{k,\lambda}}(t-\tau)^{\beta_{k,\lambda}}}{(\alpha_{k,\lambda})!(\beta_{k,\lambda})!}$ of the estimate for $|(K^{\otimes k})_{12}(s,t,\sigma,\tau)|$, the three terms $(K^{\otimes k})_{12}\otimes_{12,1}K_1$, $(K^{\otimes k})_{12}\otimes_{12,2}K_2$, $(K^{\otimes k})_{12}\otimes_{12,12}K_{12}$ contribute, respectively, the terms

$$C^{k+1}\frac{(s-\sigma)^{\alpha_{k,\lambda}+1}(t-\tau)^{\beta_{k,\lambda}}}{(\alpha_{k,\lambda}+1)!(\beta_{k,\lambda})!},\ C^{k+1}\frac{(s-\sigma)^{\alpha_{k,\lambda}}(t-\tau)^{\beta_{k,\lambda}+1}}{(\alpha_{k,\lambda})!(\beta_{k,\lambda}+1)!},$$

$$C^{k+1}\frac{(s-\sigma)^{\alpha_{k,\lambda}+1}(t-\tau)^{\beta_{k,\lambda}+1}}{(\alpha_{k,\lambda}+1)!(\beta_{k,\lambda}+1)!}$$

to the estimation of $|(K^{\otimes(k+1)})_{12}(s,t,\sigma,\tau)|$. The pairs $(\alpha_{k+1,\lambda_i},\beta_{k+1,\lambda_i})$, $i=1,...,7$, defined by specifying that they take the values $(k-1,0),(0,k-1),(k,0),(0,k),(\alpha_{k,\lambda}+1,\beta_{k,\lambda}),(\alpha_{k,\lambda},\beta_{k,\lambda}+1),(\alpha_{k,\lambda}+1,\beta_{k,\lambda}+1)$, satisfy the conditions $\alpha_{k+1,\lambda_i}+\beta_{k+1,\lambda_i}\geq k-1$, $\max(\alpha_{k+1,\lambda_i},\beta_{k+1,\lambda_i})\leq k$ provided $\alpha_{k,\lambda}+\beta_{k,\lambda}\geq k-2$, $\max(\alpha_{k,\lambda},\beta_{k,\lambda})\leq k-1$. This is the inductive step in proving (5.19) for all k; for k=1 and k=2, the validity of the wanted inequalities is verifiable by direct calculation. ///

Theorem 5.1. We set

$$Q:=\max(A,B,1),\ m_k=\left\lfloor\frac{k}{2}-1\right\rfloor \qquad \text{--- (5.21)}$$

(where $\lfloor w\rfloor$ signifies the greatest among all integers that do not exceed the real number w, in other words $\lfloor w\rfloor$ is the "floor function" of w).
Then we have

$$|(K^{\otimes k})_{12}(s,t,\sigma,\tau)|\leq\frac{(3C)^k}{[(m_k)!][k-2-m_k]!}Q^{2(k-1)} \qquad \text{--- (5.22)}$$

Proof: Among all pairs $(\alpha_{k,\lambda}, \beta_{k,\lambda})$ that satisfy $\alpha_{k,\lambda} + \beta_{k,\lambda} = m$, the maximum of $\frac{1}{(\alpha_{k,\lambda})!(\beta_{k,\lambda})!}$ is achieved at $\alpha_{k,\lambda} = \alpha_m^* := \left\lfloor \frac{m}{2} \right\rfloor$, and that maximum value is equal to $\frac{1}{[(\alpha_m^*)!][m - \alpha_m^*]!}$. (This is tantamount to the well-known result that, in the expansion of $(x+y)^m$, the maximum of the coefficients is in the middle if m is even, and it is achieved at the two terms in the middle if m is odd.) The smallest possible value of m, in the context of the estimates for $|(K^{\otimes k})_{12}(s,t,\sigma,\tau)|$, is k–2, and the corresponding maximum over m is therefore achieved at $m_k^*$. The sum of the exponents $\alpha_{k,\lambda} + \beta_{k,\lambda}$ does not exceed 2k–2, and consequently the powers $(s-\sigma)^{\alpha_{k,\lambda}}(t-\tau)^{\beta_{k,\lambda}}$ do not exceed $Q^{2(k-1)}$. The number of terms, $N_k$, in the estimate (5.17), does not exceed $3^k$. Thus the assertion of this theorem has been proved. ///

Since the estimates on $(K^{\otimes k})_1$, $(K^{\otimes k})_2$, $(K^{\otimes k})_{12}$ are all terms of absolutely convergent series, we have the following:

Theorem 5.2. Under the conditions stated before lemma 5.1, the 3 series defining the resolvent kernel R in (5.12) converge uniformly in **D**. ///

The construction of R and the representation of $z_{(N)}$ in (5.11) imply the following:

Theorem 5.3. The solution z(s, t) of (5.1) is given by

$$z(s,t) = z_0(s,t) + (R \otimes_0 z_0)(s,t) \qquad \text{--- (5.23)}$$

The resolvent kernel R satisfies

$$K \otimes R = R \otimes K = R - K \qquad \text{--- (5.24)}$$

///

We now define and analyze an equation of a type adjoint to (5.1):

$$\zeta(s,t) = \zeta_0(s,t) + \int_s^A \zeta(\sigma,t)K_1(\sigma,t,s)d\sigma + \int_t^B \zeta(s,\tau)K_2(s,\tau,t)d\tau +$$
$$+ \int_s^A \int_t^B \zeta(\sigma,\tau)K_{12}(\sigma,\tau,s,t)d\tau d\sigma$$

--- (5.25)

where the values of $\zeta(s,t)$ are n-dimensional row vectors.

We shall prove:

Theorem 5.4. For initial function $\zeta_0$ that is continuous on **D**, and under the conditions of this section on the kernel K, the solution of (5.25) is given by

$$\zeta(s,t) = \zeta_0(s,t) + \int_s^A \zeta_0(\sigma,t)R_1(\sigma,t,s)d\sigma + \int_t^B \zeta_0(s,\tau)R_2(s,\tau,t)d\tau +$$
$$+ \int_s^A \int_t^B \zeta_0(\sigma,\tau)R_{12}(\sigma,\tau,s,t)d\tau d\sigma$$

--- (5.26)

Proof: Suppose $\zeta$ is given by (5.26); we shall show that $\zeta$ solves (5.25).
For convenience, we define the convolution $\otimes_{\tilde{0}}$ by

$$(\zeta_0 \otimes_{\tilde{0}} M)(s,t) := \int_s^A \zeta_0(\sigma,t)M_1(\sigma,t,s)d\sigma + \int_t^B \zeta_0(s,\tau)M_2(s,\tau,t)d\tau +$$
$$+ \int_s^A \int_t^B \zeta_0(\sigma,\tau)M_{12}(\sigma,\tau,s,t)d\sigma d\tau$$

--- (5.27)

with the operations corresponding to the 3 terms on the right-hand side of (5.27) denoted by $\otimes_{\tilde{0},1}$, $\otimes_{\tilde{0},2}$, $\otimes_{\tilde{0},12}$, respectively, i.e.

$$(\zeta_0 \tilde{\otimes}_{0,1} M_1)(s,t) := \int_s^A \zeta_0(\sigma,t) M_1(\sigma,t,s) d\sigma;$$

$$(\zeta_0 \tilde{\otimes}_{0,2} M_2)(s,t) := \int_t^B \zeta_0(s,\tau) M_2(s,\tau,t) d\tau;$$

$$(\zeta_0 \tilde{\otimes}_{0,12} M_{12})(s,t) := \int_s^A \int_t^B \zeta_0(\sigma,\tau) M_{12}(\sigma,\tau,s,t) d\sigma d\tau$$

--- (5.28)

With $\zeta$ given by (5.26), we calculate

$$\int_s^A \zeta(\sigma,t)K_1(\sigma,t,s)d\sigma + \int_t^B \zeta(s,\tau)K_2(s,\tau,t)d\tau + \int_s^A \int_t^B \zeta(\sigma,\tau)K_{12}(\sigma,\tau,s,t)d\sigma d\tau =$$

$$= \int_s^A \zeta_0(\sigma,t)K_1(\sigma,t,s)d\sigma + \int_s^A \int_\sigma^A \zeta_0(\sigma_1,t)R_1(\sigma_1,t,\sigma)K_1(\sigma,t,s)d\sigma_1 d\sigma +$$

$$+ \int_s^A \int_t^B \zeta_0(\sigma,\tau_1)R_2(\sigma,\tau_1,t)K_1(\sigma,t,s)d\tau_1 d\sigma +$$

$$+ \int_s^A \int_\sigma^A \int_\tau^B \zeta_0(\sigma_1,\tau_1)R_{12}(\sigma_1,\tau_1,\sigma,\tau)K_1(\sigma,t,s)d\tau_1 d\sigma_1 d\sigma +$$

$$+ \int_t^B \zeta_0(s,\tau)K_2(s,\tau,t)d\tau +$$

$$+ \int_t^B \int_s^A \zeta_0(\sigma_1,\tau)R_1(\sigma_1,\tau,s)K_2(s,\tau,t)d\sigma_1 d\tau +$$

$$+ \int_t^B \int_\tau^B \zeta_0(s,\tau_1)R_2(s,\tau_1,\tau)K_2(s,\tau,t)d\tau_1 d\tau +$$

$$+ \int_t^B \int_s^A \int_\tau^B \zeta_0(\sigma_1,\tau_1)R_{12}(\sigma_1,\tau_1,s,\tau)K_2(s,\tau,t)d\tau_1 d\sigma_1 d\tau +$$

$$+ \int_s^A \int_t^B \zeta_0(\sigma,\tau)K_{12}(\sigma,\tau,s,t)d\sigma d\tau +$$

$$+ \int_s^A \int_t^B \int_\sigma^A \zeta_0(\sigma_1,\tau)R_1(\sigma_1,\tau,\sigma)K_{12}(\sigma,\tau,s,t)d\sigma_1 d\tau d\sigma +$$

$$+ \int_s^A \int_t^B \int_\tau^B \zeta_0(\sigma_1,\tau_1)R_2(\sigma,\tau_1,\tau)K_{12}(\sigma,\tau,s,t)d\tau_1 d\tau d\sigma +$$

$$+ \int_s^A \int_t^B \int_\sigma^A \int_\tau^B \zeta_0(\sigma_1,\tau_1)R_{12}(\sigma_1,\tau_1,\sigma,\tau)K_{12}(\sigma,\tau,s,t)d\tau_1 d\sigma_1 d\tau d\sigma$$

--- (5.29)

It is seen from (5.29) that, if $\zeta$ is given by (5.26), then

$$(\zeta \otimes_{\tilde{0}} K)(s,t) = (\zeta_0 \otimes_{\tilde{0}} K)(s,t) + (\zeta_0 \otimes_{\tilde{0}} (R \otimes K))(s,t)$$

--- (5.30)

Indeed, those terms in (5.29) that contain $R_1, R_2, R_{12}$ make up $\zeta_0 \otimes_{\tilde{0}} (R \otimes K)$. We verify that one of those terms of (5.29) is one of the terms that make up $\zeta_0 \otimes_{\tilde{0}} (R \otimes K)$:

$$\int_s^A \int_t^B \int_\sigma^A \zeta_0(\sigma_1,\tau) R_1(\sigma_1,\tau,\sigma) K_{12}(\sigma,\tau,s,t) d\sigma_1 d\tau d\sigma =$$

$$= \int_s^A \int_t^B \int_{\sigma_1}^\sigma \zeta_0(\sigma_1,\tau) R_1(\sigma_1,\tau,\sigma) K_{12}(\sigma,\tau,s,t) d\sigma d\tau d\sigma_1 =$$

$$= \int_s^A \int_t^B \zeta_0(\sigma_1,\tau)(R_1 \otimes_{1,12} K_{12})(\sigma_1,\tau,s,t) d\tau d\sigma_1 =$$

$$= (\zeta_0 \otimes_{\tilde{0},12} (R_1 \otimes_{1,12} K_{12}))(s,t)$$

--- (5.31)

Of course, all the relevant terms of (5.29) can be handled similarly. On the basis of the identity $R \otimes K = R - K$, we have

$$\zeta_0 \otimes_{\tilde{0}} K + \zeta_0 \otimes_{\tilde{0}} (R \otimes K) = \zeta_0 \otimes_{\tilde{0}} R$$

--- (5.32)

therefore, if $\zeta$ is given by (5.26), we have

$$\zeta_0 + \zeta \otimes_{\tilde{0}} K = \zeta_0 + \zeta_0 \otimes_{\tilde{0}} K + \zeta_0 \otimes_{\tilde{0}} (R \otimes K) = \zeta_0 + \zeta_0 \otimes_{\tilde{0}} R = \zeta$$

--- (5.33)

which is another way of writing (5.25). ///

6. Representation of the co-state.

For the simplified model (2.11), we define the co-state $\psi := [\psi_0 \ \psi_1(.) \ \psi_2(.) \ \psi_{12}(.,.)]$ by

$$\psi_0 := \nabla_y F_0(y(A,B),...);$$

$$\psi_1(s) := \nabla_y F_1(s,...) + \psi_0 \nabla_y f_1(A,B,s,...) + \int_s^A \{\nabla_y F_1(\sigma,...) + \psi_0 \nabla_y f_1(A,B,\sigma,...)\} R_1(\sigma,B,s) d\sigma;$$

$$\psi_2(t) := \nabla_y F_2(t,...) + \psi_0 \nabla_y f_2(A,B,t,...) + \int_t^B \{\nabla_y F_2(\tau,...) + \psi_0 \nabla_y f_2(A,B,\tau,...\} R_1(A,\tau,t) d\tau;$$

$$\psi_{12}(s,t) := \nabla_y F_{12}(s,t,...) + \psi_1(s) \nabla_y f_2(s,B,t,...) + \psi_2(t) \nabla_y f_1(A,t,s,...) +$$

$$+ \psi_0 \nabla_y f_{12}(A,B,s,t,...) + \int_s^A \psi_1(\sigma) \nabla_y f_{12}(\sigma,B,s,t,...) d\sigma + \int_t^B \psi_2(\tau) \nabla_y f_{12}(A,\tau,s,t,...) d\tau +$$

$$+ \int_s^A \{\nabla_y F_{12}(\sigma,t,...) + \psi_1(\sigma) \nabla_y f_2(\sigma,B,t,...) + \psi_2(t) \nabla_y f_1(A,t,\sigma,...) +$$

$$+ \int_t^B \{\nabla_y F_{12}(s,\tau,...) + \psi_1(s) \nabla_y f_2(s,B,\tau,...) + \psi_2(\tau) \nabla_y f_1(A,\tau,s,...) +$$

$$+ \psi_0 \nabla_y f_{12}(A,B,\sigma,\tau,...)\} R_2(s,\tau,t) d\tau + \psi_0 \nabla_y f_{12}(A,B,\sigma,t)\} R_1(\sigma,t,s) d\sigma +$$

$$+ \int_s^A \int_t^B \{\nabla_y F_{12}(\sigma,\tau,...) + \psi_1(\sigma) \nabla_y f_2(\sigma,B,\tau,...) + \psi_2(\tau) \nabla_y f_1(A,\tau,\sigma,...) +$$

$$+ \psi_0 \nabla_y f_{12}(A,B,\sigma,\tau)\} R_{12}(\sigma,\tau,s,t,...) d\tau d\sigma +$$

$$+ \int_s^A \int_\sigma^A \psi_1(\sigma_1) \nabla_y f_{12}(\sigma_1,B,\sigma,t,...) R_1(\sigma,t,s) d\sigma_1 d\sigma +$$

$$+ \int_s^A \int_t^B \{\psi_2(\tau) \nabla_y f_{12}(A,\tau,\sigma,t,...) R_1(\sigma,t,s) +$$

$$+ \psi_1(\sigma) \nabla_y f_{12}(\sigma,B,s,\tau,...) R_2(s,\tau,t)\} d\tau d\sigma +$$

$$+\int_t^B \int_\tau^B \psi_2(\tau_1)\nabla_y f_{12}(A,\tau_1,s,\tau,...)R_2(s,\tau,t)d\tau_1 d\tau +$$
$$+\int_s^A \int_t^B \left\{ \int_\sigma^A \psi_1(\sigma_1)\nabla_y f_{12}(\sigma_1,B,\sigma,\tau,...)d\sigma_1 + \right.$$
$$\left. \int_\tau^B \psi_2(\tau_1)\nabla_y f_{12}(A,\tau_1,\sigma,\tau,...)d\tau_1 \right\} R_{12}(\sigma,\tau,s,t)d\tau d\sigma$$

--- (6.1)

where $R=(R_1,R_2,R_{12})$ is the resolvent kernel associated with

$$K(s,t,\sigma,\tau)=(K_1(s,t,\sigma),K_2(s,t,\tau),K_{12}(s,t,\sigma,\tau))\equiv$$
$$\equiv(\nabla_y f_1(s,t,\sigma,...),\nabla_y f_2(s,t,\tau,...),\nabla_y f_{12}(s,t,\sigma,\tau,...))$$

--- (6.2)

The Hamiltonians $h_0,h_1,h_2,h_{12}$ are given by (3.15), but now with $\psi_0,\psi_1,\psi_2,\psi_{12}$ defined as in (6.1). It now follows that, when $(\psi_0,\psi_1,\psi_2,\psi_{12})$ are defined as in (6.1), then they satisfy the system of Hamiltonian equations (3.14), on the basis of the results on duality of linear Goursat-Volterra problems that were proved in section 5.

The variation $\delta y(s,t)$, obtained in (4.3), can be expressed, in view of the results of section 5, as

$$\delta y(s,t)=\int_0^s \nabla_u f_1(s,t,\sigma,...)\delta u(\sigma,t)d\sigma+\int_0^t \nabla_u f_2(s,t,\tau,...)\delta u(s,\tau)d\tau+$$

$$+\int_0^s\int_0^t \nabla_u f_{12}(s,t,\sigma,\tau,...)\delta u(\sigma,\tau)d\tau d\sigma+$$

$$+\int_0^s\int_0^\sigma R_1(s,t,\sigma)\nabla_u f_1(\sigma,t,\sigma_1,...)\delta u(\sigma_1,t)d\sigma_1 d\sigma+$$

$$+\int_0^s\int_0^t R_1(s,t,\sigma)\nabla_u f_2(\sigma,t,\tau,...)\delta u(\sigma,\tau)d\tau d\sigma+$$

$$+\int_0^s\int_0^\sigma\int_0^t R_1(s,t,\sigma)\nabla_u f_{12}(\sigma,t,\sigma_1,\tau,...)\delta u(\sigma_1,\tau)d\tau d\sigma_1 d\sigma+$$

$$+\int_0^t\int_0^s R_2(s,t,\tau)\nabla_u f_1(s,\tau,\sigma,...)\delta u(\sigma,\tau)d\sigma d\tau+$$

$$+\int_0^t\int_0^\tau R_2(s,t,\tau)\nabla_u f_2(s,\tau,\tau_1,...)\delta u(s,\tau_1)d\tau_1 d\tau+$$

$$+\int_0^t\int_0^s\int_0^\tau R_2(s,t,\tau)\nabla_u f_{12}(s,\tau,\sigma,\tau_1)\delta u(\sigma,\tau_1)d\tau_1 d\sigma d\tau+$$

$$+\int_0^s\int_0^t\int_0^\sigma R_{12}(s,t,\sigma,\tau)\nabla_u f_1(\sigma,\tau,\sigma_1,...)\delta u(\sigma_1,\tau)d\sigma_1 d\tau d\sigma+$$

$$+\int_0^s\int_0^t\int_0^\tau R_{12}(s,t,\sigma,\tau)\nabla_u f_2(\sigma,\tau,\tau_1,...)\delta u(\sigma,\tau_1)d\tau_1 d\tau d\sigma+$$

$$+\int_0^s\int_0^t\int_0^\sigma\int_0^\tau R_{12}(s,t,\sigma,\tau)\nabla_u f_{12}(\sigma,\tau,\sigma_1,\tau_1,...)\delta u(\sigma_1,\tau_1)d\tau_1 d\sigma_1 d\tau d\sigma$$

--- (6.3)

In order to simplify the presentation in this part of the paper, we write down the terms for the representation of the variation of the state.

Terms in the expression for $\delta y(s,t)$ :

$$\nabla_u f_1(s,t,\sigma,\ldots)\delta u(\sigma,t) \qquad \{\sigma\}$$
$$\nabla_u f_2(s,t,\tau,\ldots)\delta u(s,\tau) \qquad \{\tau\}$$
$$\nabla_u f_{12}(s,t,\sigma,\tau,\ldots)\delta u(\sigma,\tau) \qquad \{\sigma,\tau\}$$

$$R_1(s,t,\sigma)\begin{cases}\nabla_u f_1(\sigma,t,\sigma_1,\ldots)\delta u(\sigma_1,t) & \{\sigma_1,\sigma\}\\ \nabla_u f_2(\sigma,t,\tau,\ldots)\delta u(\sigma,\tau) & \{\tau,\sigma\}\\ \nabla_u f_{12}(\sigma,t,\sigma_1,\tau,\ldots)\delta u(\sigma_1,\tau) & \{\tau,\sigma_1,\sigma\}\end{cases}$$

$$R_2(s,t,\tau)\begin{cases}\nabla_u f_1(s,\tau,\sigma,\ldots)\delta u(\sigma,\tau) & \{\sigma,\tau\}\\ \nabla_u f_2(s,\tau,\tau_1,\ldots)\delta u(s,\tau_1) & \{\tau_1,\tau\}\\ \nabla_u f_{12}(s,\tau,\sigma,\tau_1,\ldots)\delta u(\sigma,\tau_1) & \{\tau_1,\sigma,\tau\}\end{cases}$$

$$R_{12}(s,t,\sigma,\tau)\begin{cases}\nabla_u f_1(\sigma,\tau,\sigma_1,\ldots)\delta u(\sigma_1,\tau) & \{\sigma_1,\tau,\sigma\}\\ \nabla_u f_2(\sigma,\tau,\tau_1,\ldots)\delta u(\sigma,\tau_1) & \{\tau_1,\tau,\sigma\}\\ \nabla_u f_{12}(\sigma,\tau,\sigma_1,\tau_1,\ldots)\delta u(\sigma_1,\tau_1) & \{\tau_1,\sigma_1,\tau,\sigma\}\end{cases}$$

Terms in the expression for $\delta y(s,B)$:

$$\nabla_u f_1(s,B,\sigma,\ldots)\delta u(\sigma,B) \qquad \{\sigma\}$$
$$\nabla_u f_2(s,B,\tau,\ldots)\delta u(s,\tau) \qquad \{\tau\}$$
$$\nabla_u f_{12}(s,B,\sigma,\tau,\ldots)\delta u(\sigma,\tau) \qquad \{\sigma,\tau\}$$

$$R_1(s,B,\sigma)\begin{cases}\nabla_u f_1(\sigma,B,\sigma_1,\ldots)\delta u(\sigma_1,B) & \{\sigma_1,\sigma\}\\ \nabla_u f_2(\sigma,B,\tau,\ldots)\delta u(\sigma,\tau) & \{\tau,\sigma\}\\ \nabla_u f_{12}(\sigma,B,\sigma_1,\tau,\ldots)\delta u(\sigma_1,\tau) & \{\tau,\sigma_1,\sigma\}\end{cases}$$

$$R_2(s,B,\tau)\begin{cases}\nabla_u f_1(s,\tau,\sigma,\ldots)\delta u(\sigma,\tau) & \{\sigma,\tau\}\\ \nabla_u f_2(s,\tau,\tau_1,\ldots)\delta u(s,\tau_1) & \{\tau_1,\tau\}\\ \nabla_u f_{12}(s,\tau,\sigma,\tau_1,\ldots)\delta u(\sigma,\tau_1) & \{\tau_1,\sigma,\tau\}\end{cases}$$

$$R_{12}(s,B,\sigma,\tau)\begin{cases}\nabla_u f_1(\sigma,\tau,\sigma_1,\ldots)\delta u(\sigma_1,\tau) & \{\sigma_1,\tau,\sigma\}\\ \nabla_u f_2(\sigma,\tau,\tau_1,\ldots)\delta u(\sigma,\tau_1) & \{\tau_1,\tau,\sigma\}\\ \nabla_u f_{12}(\sigma,\tau,\sigma_1,\tau_1,\ldots)\delta u(\sigma_1,\tau_1) & \{\tau_1,\sigma_1,\tau,\sigma\}\end{cases}$$

Terms in the expression for $\delta y(A,t)$:

$$
\begin{array}{lc}
\nabla_u f_1(A,t,\sigma,...)\delta u(\sigma,t) & \{\sigma\} \\
\nabla_u f_2(A,t,\tau,...)\delta u(A,\tau) & \{\tau\} \\
\nabla_u f_{12}(A,t,\sigma,\tau,...)\delta u(\sigma,\tau) & \{\sigma,\tau\}
\end{array}
$$

$$
R_1(A,t,\sigma)\begin{cases}
\nabla_u f_1(\sigma,t,\sigma_1,...)\delta u(\sigma_1,t) & \{\sigma_1,\sigma\} \\
\nabla_u f_2(\sigma,t,\tau,...)\delta u(\sigma,\tau) & \{\tau,\sigma\} \\
\nabla_u f_{12}(\sigma,t,\sigma_1,\tau,...)\delta u(\sigma_1,\tau) & \{\tau,\sigma_1,\sigma\}
\end{cases}
$$

$$
R_2(A,t,\tau)\begin{cases}
\nabla_u f_1(A,\tau,\sigma,...)\delta u(\sigma,\tau) & \{\sigma,\tau\} \\
\nabla_u f_2(A,\tau,\tau_1,...)\delta u(A,\tau_1) & \{\tau_1,\tau\} \\
\nabla_u f_{12}(A,\tau,\sigma,\tau_1,...)\delta u(\sigma,\tau_1) & \{\tau_1,\sigma,\tau\}
\end{cases}
$$

$$
R_{12}(A,t,\sigma,\tau)\begin{cases}
\nabla_u f_1(\sigma,\tau,\sigma_1,...)\delta u(\sigma_1,\tau) & \{\sigma_1,\tau,\sigma\} \\
\nabla_u f_2(\sigma,\tau,\tau_1,...)\delta u(\sigma,\tau_1) & \{\tau_1,\tau,\sigma\} \\
\nabla_u f_{12}(\sigma,\tau,\sigma_1,\tau_1,...)\delta u(\sigma_1,\tau_1) & \{\tau_1,\sigma_1,\tau,\sigma\}
\end{cases}
$$

Terms in the expression for $\delta y(A,B)$:

$$
\begin{array}{lc}
\nabla_u f_1(A,B,\sigma,...)\delta u(\sigma,B) & \{\sigma\} \\
\nabla_u f_2(A,B,\tau,...)\delta u(A,\tau) & \{\tau\} \\
\nabla_u f_{12}(A,B,\sigma,\tau,...)\delta u(\sigma,\tau) & \{\sigma,\tau\}
\end{array}
$$

$$R_1(A,B,\sigma)\begin{cases}\nabla_u f_1(\sigma,B,\sigma_1,...)\delta u(\sigma_1,B) & \{\sigma_1,\sigma\}\\ \nabla_u f_2(\sigma,t,\tau,...)\delta u(\sigma,\tau) & \{\tau,\sigma\}\\ \nabla_u f_{12}(\sigma,t,\sigma_1,\tau,...)\delta u(\sigma_1,\tau) & \{\tau,\sigma_1,\sigma\}\end{cases}$$

$$R_2(A,B,\tau)\begin{cases}\nabla_u f_1(A,\tau,\sigma,...)\delta u(\sigma,\tau) & \{\sigma,\tau\}\\ \nabla_u f_2(A,\tau,\tau_1,...)\delta u(A,\tau_1) & \{\tau_1,\tau\}\\ \nabla_u f_{12}(A,\tau,\sigma,\tau_1,...)\delta u(\sigma,\tau_1) & \{\tau_1,\sigma,\tau\}\end{cases}$$

$$R_{12}(A,B,\sigma,\tau)\begin{cases}\nabla_u f_1(\sigma,\tau,\sigma_1,...)\delta u(\sigma_1,\tau) & \{\sigma_1,\tau,\sigma\}\\ \nabla_u f_2(\sigma,\tau,\tau_1,...)\delta u(\sigma,\tau_1) & \{\tau_1,\tau,\sigma\}\\ \nabla_u f_{12}(\sigma,\tau,\sigma_1,\tau_1,...)\delta u(\sigma_1,\tau_1) & \{\tau_1,\sigma_1,\tau,\sigma\}\end{cases}$$

Next, we introduce the following notation:

$$u_{[0]} := col(u_{1,1}, u_{2,2}, u_{12,12});\ u_{[1]} := col(u_1, u_{2,2}, u_{12,2});$$
$$u_{[2]} := col(u_{1,1}, u_2, u_{12,2});\ u_{[12]} := col(u_1, u_2, u_{12})$$

--- (6.4)

We claim that the variation $\delta J$ can be represented as

$$\delta J = \int_0^A \int_0^B \nabla_{u_{[12]}} h_{12}(s,t,y(s,t),u_1(s),u_2(t),u_{12}(s,t),\psi_0,\psi_1(.),\psi_2(.),\psi_{12}(.,.))\cdot$$
$$\cdot\,\delta u_{[12]}(s,t)dtds +$$
$$+\int_0^A \nabla_{u_{[1]}} h_1(s,y(s,B),u_1(s),u_2(B),u_{12}(s,B),\psi_0,\psi_1(.))\delta u_{[1]}(s)ds +$$
$$+\int_0^B \nabla_{u_{[2]}} h_2(t,y(A,t),u_1(A),u_2(t),u_{12}(A,t),\psi_0,\psi_2(.))\delta u_{[2]}(t)dt +$$
$$+\nabla_{u_{[0]}} h_0(y(A,B),u_1(A),u_2(B),u_{12}(A,B),\psi_0)\delta u_{[0]}$$

--- (6.5)

The proof of (6.5) relies on taking inventory of all the terms that appear in the expression for $\delta J$ in (4.4), using the expressions for $\delta y$ in tables 6.1 through 6.4 above, and comparing with the analogous inventory for the terms on the right-hand side of (6.5). Because of the large number of terms involved, we shall illustrate this comparison for only a few of the relevant terms.

The term $\int_0^A \int_0^B (\nabla_y F_{12}(s,t,y(s,t),u(s,t))\delta y(s,t)dtds$ out of (4.4) involves all the terms for $\delta y(s,t)$. Out of those terms, we look at these 3 terms:

$$R_{12}(s,t,\sigma,\tau) \begin{cases} \nabla_u f_1(\sigma,\tau,\sigma_1,...)\delta u(\sigma_1,\tau) & \{\sigma_1,\tau,\sigma\} \\ \nabla_u f_2(\sigma,\tau,\tau_1,...)\delta u(\sigma,\tau_1) & \{\tau_1,\tau,\sigma\} \\ \nabla_u f_{12}(\sigma,\tau,\sigma_1,\tau_1,...)\delta u(\sigma_1,\tau_1) & \{\tau_1,\sigma_1,\tau,\sigma\} \end{cases}$$

.

The term $R_{12}(s,t,\sigma,\tau)\ \nabla_u f_{12}(\sigma,\tau,\sigma_1,\tau_1,...)\delta u(\sigma_1,\tau_1)$ , multiplied by $\nabla_y F_{12}(s,t,...)$, is matched with the term $\nabla_y F_{12}(\sigma,\tau,...)R_{12}(\sigma,\tau,s,t,...)$ , inside the integral $\int_s^A \int_t^B$ in the definition (6.1) of $\psi_{12}$, multiplied by $f_{12}$ in the definition of $h_{12}$ and then differentiated with respect to u. The term $R_{12}(s,t,\sigma,\tau)\nabla_u f_1(\sigma,\tau,\sigma_1,...)\delta u(\sigma_1,\tau)$, multiplied by $\nabla_y F_{12}(s,t,...)$, is matched with the term $\nabla_y F_{12}(\sigma,\tau,...)R_{12}(\sigma,\tau,s,t,...)$ , inside the integral $\int_s^A \int_t^B$ in the definition (6.1) of $\psi_{12}$, multiplied by $f_1$ to produce another one of the terms in the definition of $h_{12}$, and then differentiated with respect to u. The term $R_{12}(s,t,\sigma,\tau)\nabla_u f_2(\sigma,\tau,\tau_1,...)\delta u(\sigma,\tau_1)$ is analogous to the term $R_{12}(s,t,\sigma,\tau)\nabla_u f_1(\sigma,\tau,\sigma_1,...)\delta u(\sigma_1,\tau)$.

Other terms are matched in a similar way.

Therefore, we have proved:

<u>Theorem 6.1.</u> For the simplified model (2.11) and under the conditions stated in this section, the co-state $(\psi_0,\psi_1,\psi_2,\psi_{12})$ satisfies equation (3.14), and the variation of the cost functional J is given by (6.5). ///

## 7. Partial extremal principles.

The optimal control problem with state equation (2.11), and, *a fortiori*, the more general model (2.1), does not generally admit the derivation of necessary optimality conditions of the type of an extremal principle akin to the famous maximum principle of Pontryagin et al. [PBMG]. This is due to the fact that the controls that depend on a one-dimensional variable (only s or only t) appear both in single integrals and in a double integral, in the definition of the cost functional J, as well as in the variational equation (6.5); as a consequence, the integral expression (6.5) cannot be simultaneously reduced to pointwise form in the entire triple $(u_1(s), u_2(t), u_{12}(s,t))$. If we were to take, simultaneously, needle-shaped variations in all three controls, the passage to pointwise form would require division by $\delta s$ or $\delta t$ in the single integrals, and passage to the limit as $\delta s \to 0^+$ or $\delta t \to 0^+$, separately for each increment $\delta s$ or $\delta t$, but it would require division by the product $(\delta s)(\delta t)$ and passage to the limit as $(\delta s, \delta t) \to (0^+, 0^+)$ simultaneously for the Hamiltonian $h_{12}$ inside the double integral. The incompatibility of the two types of operations makes it impossible to obtain a pointwise form of the variational equation (6.5) jointly in all components of the control.
The only case in which it is possible to have an extremum principle akin to Pontryagin's maximum principle is the case in which only the control $u_{12}$ is present.
In general, we can have only a weaker statement of an extremum principle; we shall call this result <u>the partial extremum principle</u>. This result states that, if
$(u_1^*(s), u_2^*(t), u_{12}^*(s,t), u_1^*(A), u_1^*(B), u_{12}^*(s,B), u_{12}^*(A,t), u_{12}^*(A,B))$
is an optimal control, $(y^*(s,t), y^*(s,B), y^*(A,t), y^*(A,B))$ is the corresponding optimal trajectory (solution of the state equation with control
$(u_1^*(s), u_2^*(t), u_{12}^*(s,t), u_1^*(A), u_1^*(B), u_{12}^*(s,B), u_{12}^*(A,t), u_{12}^*(A,B))$), and
$(\psi_0^*, \psi_1^*(s), \psi_2^*(t), \psi_{12}^*(s,t))$ is the associated co-state, then, at points of continuity of
$(u_1^*(s), u_2^*(t), u_{12}^*(s,t), u_1^*(A), u_1^*(B), u_{12}^*(s,B), u_{12}^*(A,t), u_{12}^*(A,B))$ :

if $f_0, f_2$ are independent of $u_1$, then $u_1^*(s)$ minimizes, over the variable $u_1$, the function

$$h_1(s, y^*(s,B), u_1, u_2^*(B), u_{12}^*(s,B), \psi_0^*, \psi_1^*(.)) + \\ + \int_0^B h_{12}(s,t,y^*(s,t), u_1, u_2^*(t), u_{12}^*(s,t), \psi_0^*, \psi_1^*(.), \psi_2^*(.), \psi_{12}^*(.,.))dt$$

;

if $f_0, f_1$ are independent of $u_2$, then $u_2^*(t)$ minimizes, over the variable $u_2$, the function

$$h_2(t,y^*(A,t),u_1^*(A),u_2,u_{12}^*(A,t),\psi_0^*,\psi_2^*(.))+$$
$$+\int_0^A h_{12}(s,t,y^*(s,t).u_1^*(s),u_2,u_{12}^*(s,t),\psi_0^*,\psi_1(.)^*,\psi_2^*(.),\psi_{12}^*(.,.))ds$$

;

if $f_0, f_1, f_2$ are independent of $u_{12}$, then $u_{12}^*(s,t)$ minimizes, over the variable $u_{12}$, the function $h_{12}(s,t,y^*(s,t).u_1^*(s),u_2^*(t),u_{12},\psi_0^*,\psi_1(.)^*,\psi_2^*(.),\psi_{12}^*(.,.))$ ;

if $f_0$ is independent of $u_{12}$, then
$u_{12}^*(A,t)$ minimizes, over the variable $u_{12,1}$, the function
$h_2(t,y^*(A,t),u_1^*(A),u_2^*(t),u_{12,1},\psi_0^*,\psi_2^*(.))$ ;

if $f_0$ is independent of $u_{12}$, then $u_{12}^*(s,B)$ minimizes, over the variable $u_{12,2}$, the function $h_1(s,y^*(s,B),u_1^*(s),u_2^*(B),u_{12,2},\psi_0^*,\psi_1^*(.))$ ;

$u_1^*(A)$ minimizes, over the variable $u_{1,1}$, the function
$$h_0(y^*(A,B),u_{1,1},u_2^*(B),u_{12}^*(A,B))+$$
$$+\int_0^B F_2(t,y^*(A,t),u_{1,1},u_2^*(t),u_{12}^*(A,t))dt \quad ;$$

$u_2^*(B)$ minimizes, over the variable $u_{2,2}$, the function
$$h_0(y^*(A,B),u_1^*(A),u_{2,2},u_{12}^*(A,B))+$$
$$+\int_0^A F_1(s,y^*(s,B),u_1^*(s),u_{2,2},u_{12}^*(s,B))ds \quad ;$$

$u_{12}^*(A,B)$ minimizes, over the variable $u_{12,12}$, the function
$h_0(y^*(A,B),u_1^*(A),u_2^*(B),u_{12,12})$ .

We shall prove one of the above cases of a partial extremum principle, the proofs in the remaining cases being similar. We prove that, under the stated conditions, $u_{12}^*(s,t)$ minimizes, over the variable $u_{12}$, the function
$h_{12}(s,t,y^*(s,t).u_1^*(s),u_2^*(t),u_{12},\psi_0^*,\psi_1(.)^*,\psi_2^*(.),\psi_{12}^*(.,.))$ .

8. Proof of the extremum principle for the control $u_{12}(s,t)$.

We structure this proof so that it shows the significance of the condition that $f_0, f_1, f_2$ are independent of $u_{12}$. Thus, we start without that condition, and we proceed up to a certain point where that condition needs to be introduced. We make use of some ideas from [GK]; we denote that [GK] treats optimal control of ordinary differential equations, and there is no straightforward way to extend those methods to Goursat-Volterra systems. When all other controls $u_1(s), u_2(t), u_1(A), u_2(B), u_{12}(s,B), u_{12}(A,t), u_{12}(A,B)$ are set to optimal values and only the control $u_{12}(s,t)$, $0 < s < A$, $0 < t < B$ is to be determined, we rewrite the state dynamics, the cost functional, and the Hamiltonian equations in relevant form, i.e. we suppress the dependence on controls that have been set to optimal values, and we show only the dependence on the control $u_{12}(s,t)$. In this section, we shall also denote $u_{12}(s,t)$ by simply u(s, t). Thus we have the state dynamics

$$y(s,t) = g_0(s,t,u(s,t)) + \int_0^s g_1(s,t,\sigma,y(\sigma,t),u(\sigma,t))\,d\sigma + \\ + \int_0^t g_2(s,t,\tau,y(s,\tau),u(s,\tau))\,d\tau + \int_0^s \int_0^t g_{12}(s,t,\sigma,\tau,y(\sigma,\tau),u(\sigma,\tau))\,d\tau\,d\sigma \qquad \text{--- (8.1)}$$

The cost functional is now expressed as

$$J := G_0(y(A,B)) + \int_0^A G_1(s,y(s,B))\,ds + \int_0^B G_2(t,y(A,t))\,dt + \\ + \int_0^A \int_0^B G_{12}(s,t,y(s,t),u(s,t))\,dt\,ds \qquad \text{--- (8.2)}$$

The Hamiltonian equations are

$$\psi_{12}(s,t) = \nabla_y G_{12}(s,t,y(s,t),u(s,t)) + \psi_1(s)\nabla_y g_2(s,B,t,y(s,t),u(s,t)) +$$
$$+ \psi_2(t)\nabla_y g_1(A,t,s,y(s,t),u(s,t)) + \psi_0 \nabla_y g_{12}(A,B,s,t,y(s,t),u(s,t)) +$$
$$+ \int_s^A \left[ \psi_{12}(\sigma,t)\nabla_y g_1(\sigma,t,s,y(s,t),u(s,t)) + \psi_1(\sigma)\nabla_y g_{12}(\sigma,B,s,t,y(s,t),u(s,t)) \right] d\sigma +$$
$$+ \int_t^B \left[ \psi_{12}(s,\tau)\nabla_y g_2(s,\tau,t,y(s,t),u(s,t)) + \psi_2(\tau)\nabla_y g_{12}(A,\tau,s,t,y(s,t),u(s,t)) \right] d\tau +$$
$$+ \int_s^A \int_t^B \psi_{12}(\sigma,\tau)\nabla_y g_{12}(\sigma,\tau,s,t,y(s,t),u(s,t))\, d\tau\, d\sigma;$$

$$\psi_1(s) = \nabla_y G_1(s,y(s,B)) + \psi_0 \nabla_y g_1(A,B,s,y(s,B),u(s,B)) +$$
$$+ \int_s^A \psi_1(\sigma)\nabla_y g_1(\sigma,B,s,y(s,B),u(s,B))\, d\sigma \;;$$
$$\psi_2(t) = \nabla_y G_2(t,y(A,t)) + \psi_0 \nabla_y g_2(A,B,t,y(A,t),u(A,t)) +$$
$$+ \int_t^B \psi_2(\tau)\nabla_y g_2(A,\tau,t,y(A,t),u(A,t))\, d\tau \;;$$
$$\psi_0 = \nabla_y G_0(y(A,B))$$

--- (8.3)

The corresponding Hamiltonians are denoted, in this section, by $\eta_0, \eta_1, \eta_2, \eta_{12}$, thus

$$\eta_0 := G_0(y(A,B));$$

$$\eta_1(s,y,\psi_0,\psi_1(.)) := G_1(s,y) + \psi_0 g_1(A,B,s,y,u^*(s,B)) + \int_s^A \psi_1(\sigma) g_1(\sigma,B,s,y,u^*(s,B))\, d\sigma;$$

$$\eta_2(t,y,\psi_0,\psi_2(.)) := G_2(t,y) + \psi_0 g_2(A,B,t,y,u^*(A,t)) + \int_t^B \psi_2(\tau) g_2(\tau,B,s,y,u^*(A,t))\, d\tau;$$

$$\eta_{12}(s,t,y,\psi_0,\psi_1(.),\psi_2(.),\psi_{12}(.,.),u) :=$$
$$= G_{12}(s,t,y,u) + \psi_1(s) g_2(s,B,t,y,u) + \psi_2(t) g_1(A,t,s,y,u) + \psi_0 g_{12}(A,B,s,t,y,u) +$$
$$+ \int_s^A \left[ \psi_{12}(\sigma,t) g_1(\sigma,t,s,y,u) + \psi_1(\sigma) g_{12}(\sigma,B,s,t,y,u) \right] d\sigma +$$
$$+ \int_t^B \left[ \psi_{12}(s,\tau) g_2(s,\tau,t,y,u) + \psi_2(\tau) g_{12}(A,\tau,s,t,y,u) \right] d\tau +$$
$$+ \int_s^A \int_t^B \psi_{12}(\sigma,\tau) g_{12}(\sigma,\tau,s,t,y,u)\, d\tau\, d\sigma;$$

--- (8.4)

The proof consists of the following parts:
Part I: properties of the finite variations of the state, the cost functional, and the Hamiltonian;
Part II: construction of "thin" variations of the control;
Part III: construction of multi-exponential functions, and proof of a two-dimensional version of the Gronwall inequality;
Part IV: passing to the limit and completion of the proof of the extremum principle.

Part I.

We assume adequate differentiability of all functions involved.

We denote by $(u,y),(u^{\sim},y^{\sim})$ two pairs of state and control.

We show explicitly the state and the control in the functional J. Thus, for any continuous function z(.,.), and any piecewise contninuous function v(.,.), where z(.,.) is not necessarily solution of the state dynamics with control v(.,.), we set

$$J(z,v) := G_0(z(A,B)) + \int_0^A G_1(s,z(s,B))\, ds + \int_0^B G_2(t,z(A,t))\, dt +$$
$$+ \int_0^A \int_0^B G_{12}(s,t,z(s,t),v(s,t))\, dt\, ds$$

--- (8.5)

We conventionally designate the pair (y, u) as the reference pair, and the pair $(y^{\sim}, u^{\sim})$ as the comparison pair.
The finite variation of the functional J is denoted by $\Delta J(y,u)$ and is defined by

$$\Delta J(y,u) := J(y^{\sim}, u^{\sim}) - J(y,u)$$

--- (8.6)

The partial finite variations are defined by

$$\Delta_{y^{\sim}} J(y,u) := J(y^{\sim}, u) - J(y,u); \; \Delta_{u^{\sim}} J(y,u) := J(y, u^{\sim}) - J(y,u)$$

--- (8.7)

The same notation applies to other functions or functionals. For example, for the Hamiltonian $\eta_{12}(s,t,y,\psi_0,\psi_1(.),\psi_2(.),\psi_{12}(.,.),u)$, we have

$$\Delta\eta_{12}(s,t,y,\psi_0,\psi_1(.),\psi_2(.),\psi_{12}(.,.),u) :=$$
$$= \eta_{12}(s,t,y^{\sim},\psi_0,\psi_1(.),\psi_2(.),\psi_{12}(.,.),u^{\sim}) - \eta_{12}(s,t,y,\psi_0,\psi_1(.),\psi_2(.),\psi_{12}(.,.),u);$$
$$\Delta_{y^{\sim}}\eta_{12}(s,t,y,\psi_0,\psi_1(.),\psi_2(.),\psi_{12}(.,.),u) :=$$
$$= \eta_{12}(s,t,y^{\sim},\psi_0,\psi_1(.),\psi_2(.),\psi_{12}(.,.),u) - \eta_{12}(s,t,y,\psi_0,\psi_1(.),\psi_2(.),\psi_{12}(.,.),u);$$
$$\Delta_{u^{\sim}}\eta_{12}(s,t,y,\psi_0,\psi_1(.),\psi_2(.),\psi_{12}(.,.),u) :=$$
$$= \eta_{12}(s,t,y,\psi_0,\psi_1(.),\psi_2(.),\psi_{12}(.,.),u^{\sim}) - \eta_{12}(s,t,y,\psi_0,\psi_1(.),\psi_2(.),\psi_{12}(.,.),u)$$

--- (8.8)

We shall use the following formula for the variation of a differentiable (with respect to y) function $\varphi(s,t,y(s,t),u(s,t))$:

$$\Delta\varphi(s,t,y(s,t),u(s,t)) = \Delta_{u^{\sim}}\varphi(s,t,y(s,t),u(s,t)) +$$
$$+ [\nabla_y \varphi(s,t,y(s,t),u(s,t)) + \nabla_y \Delta_{u^{\sim}}\varphi(s,t,y(s,t),u(s,t))]\Delta y(s,t) + \mathbf{o}_{\varphi}(\| \Delta y(s,t) \|)$$

--- (8.9)

This is proved as follows:

$$\Delta\varphi(s,t,y(s,t),u(s,t)) = \varphi(s,t,y^{\sim}(s,t),u^{\sim}(s,t)) - \varphi(s,t,y(s,t),u(s,t)) =$$
$$= \varphi(s,t,y^{\sim}(s,t),u^{\sim}(s,t)) - \varphi(s,t,y^{\sim}(s,t),u(s,t)) + \varphi(s,t,y^{\sim}(s,t),u(s,t)) -$$
$$- \varphi(s,t,y(s,t),u(s,t)) =$$
$$= \Delta_{u^{\sim}}\varphi(s,t,y^{\sim}(s,t),u(s,t)) + \Delta_{y^{\sim}}\varphi(s,t,y(s,t),u(s,t)) =$$
$$= \Delta_{u^{\sim}}\varphi(s,t,y(s,t),u(s,t)) + \nabla_y\Delta_{u^{\sim}}\varphi(s,t,y(s,t),u(s,t))\Delta y(s,t) + \mathbf{o}_{\varphi,1}(\|\Delta y(s,t)\|) +$$
$$+ \nabla_y\varphi(s,t,y(s,t),u(s,t))\Delta y(s,t) + \mathbf{o}_{\varphi,2}(\|\Delta y(s,t)\|)$$
--- (8.10)

thus we obtain (8.9) with $\mathbf{o}_{\varphi}(\|\Delta y(s,t)\|) = \mathbf{o}_{\varphi,1}(\|\Delta y(s,t)\|) + \mathbf{o}_{\varphi,2}(\|\Delta y(s,t)\|)$. (The symbols $\mathbf{o}_j(\xi)$, for different expressions substituted in lieu of the subscript j, signify different quantities that are $\mathbf{o}(\xi)$ as $\xi \to 0^+$, where $\mathbf{o}$ stands for Landau's little-oh symbol.)

Now, we claim that the finite variation of the functional J can be represented as follows:

$$\Delta J(y,u) = \int_0^A \int_0^B \{\Delta_{u^{\sim}}\eta_{12}(s,t,y(s,t),\psi_0,\psi_1(.),\psi_2(.),\psi_{12}(.,.),u(s,t)) +$$
$$+ \nabla_y\Delta_{u^{\sim}}\eta_{12}(s,t,y(s,t),\psi_0,\psi_1(.),\psi_2(.),\psi_{12}(.,.),u(s,t))\Delta y(s,t) + \mathbf{o}_{J,12}(\|\Delta y(s,t)\|)\}\,dt\,ds +$$
$$+ \int_0^A \mathbf{o}_{J,1}(\|\Delta y(s,B)\|)\,ds + \int_0^B \mathbf{o}_{J,2}(\|\Delta y(A,t)\|)\,dt + \mathbf{o}_{J,0}(\|\Delta y(A,B)\|)$$
--- (8.11)

To prove (8.11), we first calculate, in accordance with (8.9),

$$\Delta J(y,u) = \nabla_y G_0(y(A,B))\Delta y(A,B) + \mathbf{o}_{J,0,1}(\|\Delta y(A,B)\|) +$$
$$+ \int_0^A \{\nabla_y G_1(s,y(s,B))\Delta y(s,B) + \mathbf{o}_{J,1,1}(\|\Delta y(s,B)\|)\}\,ds +$$
$$+ \int_0^B \{\nabla_y G_2(t,y(A,t))\Delta y(A,t) + \mathbf{o}_{J,2,1}(\|\Delta y(A,t)\|)\}\,dt +$$
$$+ \int_0^A \int_0^B \{\Delta_{u^{\sim}} G_{12}(s,t,y(s,t),u(s,t)) + [\nabla_y G_{12}(s,t,y(s,t),u(s,t)) +$$
$$+ \nabla_y\Delta_{u^{\sim}} G_{12}(s,t,y(s,t),u(s,t))]\Delta y(s,t) + \mathbf{o}_{J,12,1}(\|\Delta y(s,t)\|)\}\,dt\,ds$$

--- (8.12)

The quantities $\Delta y(s,t), \Delta y(s,B), \Delta y(A,t), \Delta y(A,B)$ are expressed, on the basis of (8.1) and (8.9), as follows:

$$\Delta y(s,t) = \Delta_{u^{\sim}} g_0(s,t,u(s,t)) + \int_0^s \{\Delta_{u^{\sim}} g_1(s,t,\sigma,y(\sigma,t),u(\sigma,t)) +$$
$$+ [\nabla_y g_1(s,t,\sigma,y(\sigma,t),u(\sigma,t)) + \nabla_y \Delta_{u^{\sim}} g_1(s,t,\sigma,y(\sigma,t),u(\sigma,t))]\Delta y(\sigma,t) +$$
$$+ \mathbf{o}_{y,12,1}(\|\Delta y(\sigma,t)\|)\}\, d\sigma +$$
$$+ \int_0^t \{\Delta_{u^{\sim}} g_2(s,t,\tau,y(s,\tau),u(s,\tau)) + [\nabla_y g_2(s,t,\tau,y(s,\tau),u(s,\tau)) +$$
$$+ \nabla_y \Delta_{u^{\sim}} g_2(s,t,\tau,y(s,\tau),u(s,\tau))]\Delta y(s,\tau) + \mathbf{o}_{y,12,2}(\|\Delta y(s,\tau)\|)\}\, d\tau +$$
$$+ \int_0^s \int_0^t \{\Delta_{u^{\sim}} g_{12}(s,t,\sigma,\tau,y(\sigma,\tau),u(\sigma,\tau)) + [\nabla_y g_{12}(s,t,\sigma,\tau,y(\sigma,\tau),u(\sigma,\tau)) +$$
$$+ \nabla_y \Delta_{u^{\sim}} g_{12}(s,t,\sigma,\tau,y(\sigma,\tau),u(\sigma,\tau))]\Delta y(\sigma,\tau) + \mathbf{o}_{y,12,12}(\|\Delta y(\sigma,\tau)\|)\}\, d\tau\, d\sigma$$

--- (8.13)

$$\Delta y(s,B) = \int_0^s \{\nabla g_1(s,B,y(\sigma,B),u^*(\sigma,B)) + \mathbf{o}_{y,1,1}(\|\Delta y(\sigma,B)\|)\} d\sigma +$$
$$+ \int_0^B \{\Delta_{u^{\sim}} g_2(s,B,\tau,y(s,\tau),u(s,\tau)) + [\nabla_y g_2(s,B,\tau,y(s,\tau),u(s,\tau)) +$$
$$+ \nabla_y \Delta_{u^{\sim}} g_2(s,B,\tau,y(s,\tau),u(s,\tau))]\Delta y(s,\tau) + \mathbf{o}_{y,1,2}(\|\Delta y(s,\tau)\|)\}\, d\tau +$$
$$+ \int_0^s \int_0^B \{\Delta_{u^{\sim}} g_{12}(s,B,\sigma,\tau,y(\sigma,\tau),u(\sigma,\tau)) + [\nabla_y g_{12}(s,B,\sigma,\tau,y(\sigma,\tau),u(\sigma,\tau)) +$$
$$+ \nabla_y \Delta_{u^{\sim}} g_{12}(s,B,\sigma,\tau,y(\sigma,\tau),u(\sigma,\tau))]\Delta y(\sigma,\tau) + \mathbf{o}_{y,1,12}(\|\Delta y(\sigma,\tau)\|)\}\, d\tau\, d\sigma$$

--- (8.14)

$$\Delta y(A,t)=\int_0^A \{\Delta_{\tilde u} g_1(A,t,\sigma,y(\sigma,t),u(\sigma,t))+[\nabla_y g_1(A,t,\sigma,y(\sigma,t),u(\sigma,t))+$$
$$+\nabla_y\Delta_{\tilde u} g_1(A,t,\sigma,y(\sigma,t),u(\sigma,t))]\Delta y(A,t)+\mathbf{o}_{y,2,1}(\|\Delta y(\sigma,t)\|)\}\,d\sigma+$$
$$+\int_0^t \{\nabla_y g_2(A,t,\tau,y(A,\tau),u^*(A,\tau))\Delta y(A,\tau)+\mathbf{o}_{y,2,2}(\|\Delta y(A,\tau)\|)\}\,d\tau+$$
$$+\int_0^A\int_0^t \{\Delta_{\tilde u} g_{12}(A,t,\sigma,\tau,y(\sigma,\tau),u(\sigma,\tau))+[\nabla_y g_{12}(A,t,\sigma,\tau,y(\sigma,\tau),u(\sigma,\tau))+$$
$$+\nabla_y\Delta_{\tilde u} g_{12}(A,t,\sigma,\tau,y(\sigma,\tau),u(\sigma,\tau))]\Delta y(\sigma,\tau)+\mathbf{o}_{y,2,12}(\|\Delta y(\sigma,\tau)\|)\}\,d\tau\, d\sigma$$
--- (8.15)

$$\Delta y(A,B)=\int_0^A \{\nabla_y g_1(A,B,\sigma,y(\sigma,B),u^*(\sigma,B))\Delta y(\sigma,B)+\mathbf{o}_{y,0,1}(\|\Delta y(\sigma,B)\|)\}\,d\sigma+$$
$$+\int_0^B \{\nabla_y g_2(A,B,\tau,y(A,\tau),u^*(A,\tau))\Delta y(A,\tau)+\mathbf{o}_{y,0,2}(\|\Delta y(A,t)\|)\}\,d\tau+$$
$$+\int_0^A\int_0^B \{\Delta_{\tilde u} g_{12}(A,B,\sigma,\tau,y(\sigma,\tau),u(\sigma,\tau))+[\nabla_y g_{12}(A,B,\sigma,\tau,y(\sigma,\tau),u(\sigma,\tau))+$$
$$+\nabla_y\Delta_{\tilde u} g_{12}(A,B,\sigma,\tau,y(\sigma,\tau),u(\sigma,\tau))]\Delta y(\sigma,\tau)+\mathbf{o}_{y,0,12}(\|\Delta y(\sigma,\tau)\|)\}\,d\tau\, d\sigma$$
--- (8.16)

On the basis of formulae (8.13) through (8.16), we evaluate the expression

$$E:=\int_0^A\int_0^B \psi_{12}(s,t)\Delta y(s,t)\,dt\,ds+\int_0^A \psi_1(s)\Delta y(s,B)\,ds+\int_0^B \psi_2(t)\Delta y(A,t)\,dt+\psi_0\Delta y(A,B)$$
--- (8.17)

The quantity E plays a role similar to that of $\int_{t_0}^{t_1} \psi(t)\Delta x(t)\,dt$ in the theory of Pontryagin's maximum principle for controlled ordinary differential equations, as developed in [GK]; E is also related, but not exactly analogous, to the energy of a mechanical system that is described by the equations of Hamiltonian dynamics. The properties of E are established next.

We have

$$\int_0^A \int_0^B \psi_{12}(s,t)\Delta y(s,t)\,dt\,ds = \int_0^A \int_0^B \psi_{12}(s,t)\Delta_{\tilde{u}} g_0(s,t,u(s,t))\,ds\,dt +$$
$$+\int_0^A \int_0^B \int_0^s \psi_{12}(s,t)\{\Delta_{\tilde{u}} g_1(s,t,\sigma,y(\sigma,t),u(\sigma,t)) +$$
$$+[\nabla_y g_1(s,t,\sigma,y(\sigma,t),u(\sigma,t)) + \nabla_y \Delta_{\tilde{u}} g_1(s,t,\sigma,y(\sigma,t),u(\sigma,t))]\Delta y(\sigma,t) +$$
$$+\mathbf{o}_{y,12,1}(\| \Delta y(\sigma,t) \|)\}\,d\sigma\,dt\,ds +$$
$$+\int_0^A \int_0^B \int_0^t \psi_{12}(s,t)\{\Delta_{\tilde{u}} g_2(s,t,\tau,y(s,\tau),u(s,\tau)) + [\nabla_y g_2(s,t,\tau,y(s,\tau),u(s,\tau)) +$$
$$+\nabla_y \Delta_{\tilde{u}} g_2(s,t,\tau,y(s,\tau),u(s,\tau))]\Delta y(s,\tau) + \mathbf{o}_{y,12,2}(\| \Delta y(s,\tau) \|)\}\,d\tau\,dt\,ds +$$
$$+\int_0^A \int_0^B \int_0^s \int_0^t \psi_{12}(s,t)\{\Delta_{\tilde{u}} g_{12}(s,t,\sigma,\tau,y(\sigma,\tau),u(\sigma,\tau)) + [\nabla_y g_{12}(s,t,\sigma,\tau,y(\sigma,\tau),u(\sigma,\tau)) +$$
$$+\nabla_y \Delta_{\tilde{u}} g_{12}(s,t,\sigma,\tau,y(\sigma,\tau),u(\sigma,\tau))]\Delta y(\sigma,\tau) + \mathbf{o}_{y,12,12}(\| \Delta y(\sigma,\tau) \|)\}\,d\tau\,d\sigma\,dt\,ds$$

--- (8.18)

After some changes in the order of integration, we get

$$\int_0^A \int_0^B \psi_{12}(s,t)\Delta y(s,t)\,dt\,ds = \int_0^A \int_0^B \psi_{12}(s,t)\Delta_{\tilde{u}} g_0(s,t,u(s,t))\,ds\,dt +$$
$$+\int_0^A \int_0^B \int_s^A \psi_{12}(\sigma,t)\{\Delta_{\tilde{u}} g_1(\sigma,t,s,y(s,t),u(s,t)) +$$
$$+[\nabla_y g_1(\sigma,t,s,y(s,t),u(s,t)) + \nabla_y \Delta_{\tilde{u}} g_1(\sigma,t,s,y(s,t),u(s,t))]\Delta y(s,t) +$$
$$+\mathbf{o}_{y,12,1}(\| \Delta y(s,t) \|)\}\,d\sigma\,dt\,ds +$$
$$+\int_0^A \int_0^B \int_t^B \psi_{12}(s,\tau)\{\Delta_{\tilde{u}} g_2(s,\tau,t,y(s,t),u(s,t)) + [\nabla_y g_2(s,\tau,t,y(s,t),u(s,t)) +$$
$$+\nabla_y \Delta_{\tilde{u}} g_2(s,\tau,t,y(s,t),u(s,t))]\Delta y(s,t) + \mathbf{o}_{y,12,2}(\| \Delta y(s,t) \|)\}\,d\tau\,dt\,ds +$$
$$+\int_0^A \int_0^B \int_s^A \int_t^B \psi_{12}(\sigma,\tau)\{\Delta_{\tilde{u}} g_{12}(\sigma,\tau,s,t,y(s,t),u(s,t)) + [\nabla_y g_{12}(\sigma,\tau,s,t,y(s,t),u(s,t)) +$$
$$+\nabla_y \Delta_{\tilde{u}} g_{12}(\sigma,\tau,s,t,y(s,t),u(s,t))]\Delta y(s,t) + \mathbf{o}_{y,12,12}(\| \Delta y(s,t) \|)\}\,d\tau\,d\sigma\,dt\,ds$$

--- (8.19)

Similarly, we find

$$\int_0^A \psi_1(s)\Delta y(s,B)\,ds = \int_0^A \int_0^s \psi_1(s)\{\nabla g_1(s,B,y(\sigma,B),u^*(\sigma,B)) + \mathbf{o}_{y,1,1}(\|\Delta y(\sigma,B)\|)\}d\sigma\,ds +$$

$$+\int_0^A \int_0^B \psi_1(s)\{\Delta_{\tilde{u}} g_2(s,B,\tau,y(s,\tau),u(s,\tau)) + [\nabla_y g_2(s,B,\tau,y(s,\tau),u(s,\tau)) +$$

$$+\nabla_y \Delta_{\tilde{u}} g_2(s,B,\tau,y(s,\tau),u(s,\tau))]\Delta y(s,\tau) + \mathbf{o}_{y,1,2}(\|\Delta y(s,\tau)\|)\}\,d\tau\,ds +$$

$$+\int_0^A \int_0^s \int_0^B \psi_1(s)\{\Delta_{\tilde{u}} g_{12}(s,B,\sigma,\tau,y(\sigma,\tau),u(\sigma,\tau)) + [\nabla_y g_{12}(s,B,\sigma,\tau,y(\sigma,\tau),u(\sigma,\tau)) +$$

$$+\nabla_y \Delta_{\tilde{u}} g_{12}(s,B,\sigma,\tau,y(\sigma,\tau),u(\sigma,\tau))]\Delta y(\sigma,\tau) + \mathbf{o}_{y,1,12}(\|\Delta y(\sigma,\tau)\|)\}\,d\tau\,d\sigma\,ds =$$

$$= \int_0^A \int_s^A \psi_1(\sigma)\{\nabla g_1(\sigma,B,y(s,B),u^*(s,B)) + \mathbf{o}_{y,1,1}(\|\Delta y(s,B)\|)\}d\sigma\,ds +$$

$$+\int_0^A \int_0^B \psi_1(s)\{\Delta_{\tilde{u}} g_2(s,B,t,y(s,t),u(s,t)) + [\nabla_y g_2(s,B,t,y(s,t),u(s,t)) +$$

$$+\nabla_y \Delta_{\tilde{u}} g_2(s,B,t,y(s,t),u(s,t))]\Delta y(s,t) + \mathbf{o}_{y,1,2}(\|\Delta y(s,t)\|)\}\,dt\,ds +$$

$$+\int_0^A \int_0^B \int_s^A \psi_1(\sigma)\{\Delta_{\tilde{u}} g_{12}(\sigma,B,s,t,y(s,t),u(s,t)) + [\nabla_y g_{12}(\sigma,B,s,t,y(s,t),u(s,t)) +$$

$$+\nabla_y \Delta_{\tilde{u}} g_{12}(\sigma,B,s,t,y(s,t),u(s,t))]\Delta y(s,t) + \mathbf{o}_{y,1,12}(\|\Delta y(s,t)\|)\}\,d\sigma\,dt\,ds$$

--- (8.20)

$$\int_0^B \psi_2(t)\Delta y(A,t)\,dt = \int_0^B \int_0^A \psi_2(t)\{\Delta_{\tilde{u}} g_1(A,t,\sigma,y(\sigma,t),u(\sigma,t)) +$$
$$+[\nabla_y g_1(A,t,\sigma,y(\sigma,t),u(\sigma,t)) + \nabla_y \Delta_{\tilde{u}} g_1(A,t,\sigma,y(\sigma,t),u(\sigma,t))]\Delta y(\sigma,t) +$$
$$\mathbf{o}_{y,2,1}(\|\Delta y(\sigma,t)\|)\}\,d\sigma\,dt +$$
$$+\int_0^B \int_0^t \psi_2(t)\{\nabla_y g_2(A,t,\tau,y(A,\tau),u^*(A,\tau))\Delta y(A,\tau) + \mathbf{o}_{y,2,2}(\|\Delta y(A,\tau)\|)\}\,d\tau\,dt +$$
$$+\int_0^B \int_0^A \int_0^t \psi_2(t)\{\Delta_{\tilde{u}} g_{12}(A,t,\sigma,\tau,y(\sigma,\tau),u(\sigma,\tau)) + [\nabla_y g_{12}(A,t,\sigma,\tau,y(\sigma,\tau),u(\sigma,\tau)) +$$
$$+\nabla_y \Delta_{\tilde{u}} g_{12}(A,t,\sigma,\tau,y(\sigma,\tau),u(\sigma,\tau))]\Delta y(\sigma,\tau) + \mathbf{o}_{y,2,12}(\|\Delta y(\sigma,\tau)\|)\}\,d\tau\,d\sigma\,dt =$$
$$= \int_0^A \int_0^B \psi_2(t)\{\Delta_{\tilde{u}} g_1(A,t,s,y(s,t),u(s,t)) +$$
$$+[\nabla_y g_1(A,t,s,y(s,t),u(s,t)) + \nabla_y \Delta_{\tilde{u}} g_1(A,t,s,y(s,t),u(s,t))]\Delta y(s,t) +$$
$$\mathbf{o}_{y,2,1}(\|\Delta y(\sigma,t)\|)\}\,dt\,ds +$$
$$+\int_0^B \int_t^B \psi_2(\tau)\{\nabla_y g_2(A,\tau,t,y(A,t),u^*(A,t))\Delta y(A,t) + \mathbf{o}_{y,2,2}(\|\Delta y(A,t)\|)\}\,d\tau\,dt +$$
$$+\int_0^B \int_0^A \int_t^B \psi_2(\tau)\{\Delta_{\tilde{u}} g_{12}(A,\tau,s,t,y(s,t),u(s,t)) + [\nabla_y g_{12}(A,\tau,s,t,y(s,t),u(s,t)) +$$
$$+\nabla_y \Delta_{\tilde{u}} g_{12}(A,\tau,s,t,y(s,t),u(s,t))]\Delta y(s,t) + \mathbf{o}_{y,2,12}(\|\Delta y(s,t)\|)\}\,d\tau\,ds\,dt$$

--- (8.21)

$$\psi_0 \Delta y(A,B) = \int_0^A \psi_0\{\nabla_y g_1(A,B,s,y(s,B),u^*(s,B))\Delta y(s,B) + \mathbf{o}_{y,0,1}(\|\Delta y(s,B)\|)\}\,ds +$$
$$+\int_0^B \psi_0\{\nabla_y g_2(A,B,s,y(A,t),u^*(A,t))\Delta y(A,t) + \mathbf{o}_{y,0,2}(\|\Delta y(A,t)\|)\}\,dt +$$
$$+\int_0^A \int_0^B \psi_0\{\Delta_{\tilde{u}} g_{12}(A,B,s,t,y(s,t),u(s,t)) + [\nabla_y g_{12}(A,B,s,t,y(s,t),u(s,t)) +$$
$$+\nabla_y \Delta_{\tilde{u}} g_{12}(A,B,s,t,y(s,t),u(s,t))]\Delta y(s,t) + \mathbf{o}_{y,0,12}(\|\Delta y(s,t)\|)\}\,dt\,ds$$

--- (8.22)

We now use the expression (8.12) for $\Delta J$ and the Hamiltonian equations (8.3) and (8.4). Then we can verify that

$$E + \Delta J = E + \int_0^A \int_0^B \{\Delta_{u^\sim} \eta_{12}(s,t,y(s,t),\psi_0,\psi_1(.),\psi_2(.),\psi_{12}(.,.),u(s,t)) +$$

$$+ \nabla_y \Delta_{u^\sim} \eta_{12}(s,t,y(s,t),\psi_0,\psi_1(.),\psi_2(.),\psi_{12}(.,.),u(s,t))\Delta y(s,t) + \mathbf{o}_{J,12}(\| \Delta y(s,t) \|)\}\,dt\,ds +$$

$$+ \int_0^A \mathbf{o}_{J,1}(\| \Delta y(s,B) \|)\,ds + \int_0^B \mathbf{o}_{J,2}(\| \Delta y(A,t) \|)\,dt + \mathbf{o}_{J,0}(\| \Delta y(A,B) \|)$$

--- (8.23)

which is tantamount to (8.11).

Indeed,

$$E + \Delta J = \int_0^A \int_0^B \psi_{12}(s,t)\Delta_{u^\sim} g_0(s,t,u(s,t))\,ds\,dt +$$

$$= \int_0^A \int_0^B \int_s^A \psi_{12}(\sigma,t)\{\Delta_{u^\sim} g_1(\sigma,t,s,y(s,t),u(s,t)) +$$

$$+ [\nabla_y g_1(\sigma,t,s,y(s,t),u(s,t)) + \nabla_y \Delta_{u^\sim} g_1(\sigma,t,s,y(s,t),u(s,t))]\Delta y(s,t) +$$

$$+ \mathbf{o}_{y,12,1}(\| \Delta y(s,t) \|)\}\,d\sigma\,dt\,ds +$$

$$+ \int_0^A \int_0^B \int_t^B \psi_{12}(s,\tau)\{\Delta_{u^\sim} g_2(s,\tau,t,y(s,t),u(s,t)) + [\nabla_y g_2(s,\tau,t,y(s,t),u(s,t)) +$$

$$+ \nabla_y \Delta_{u^\sim} g_2(s,\tau,t,y(s,t),u(s,t))]\Delta y(s,t) + \mathbf{o}_{y,12,2}(\| \Delta y(s,t) \|)\}\,d\tau\,dt\,ds +$$

$$+ \int_0^A \int_0^B \int_s^A \int_t^B \psi_{12}(\sigma,\tau)\{\Delta_{u^\sim} g_{12}(\sigma,\tau,s,t,y(s,t),u(s,t)) + [\nabla_y g_{12}(\sigma,\tau,s,t,y(s,t),u(s,t)) +$$

$$+ \nabla_y \Delta_{u^\sim} g_{12}(\sigma,\tau,s,t,y(s,t),u(s,t))]\Delta y(s,t) + \mathbf{o}_{y,12,12}(\| \Delta y(s,t) \|)\}\,d\tau\,d\sigma\,dt\,ds +$$

$$+ \int_0^A \int_s^A \psi_1(\sigma)\{\nabla g_1(\sigma,B,y(s,B),u^*(s,B)) + \mathbf{o}_{y,1,1}(\| \Delta y(s,B) \|)\}d\sigma\,ds +$$

$$+ \int_0^A \int_0^B \psi_1(s)\{\Delta_{u^\sim} g_2(s,B,t,y(s,t),u(s,t)) + [\nabla_y g_2(s,B,t,y(s,t),u(s,t)) +$$

$$+ \nabla_y \Delta_{u^\sim} g_2(s,B,t,y(s,t),u(s,t))]\Delta y(s,t) + \mathbf{o}_{y,1,2}(\| \Delta y(s,t) \|)\}\,dt\,ds +$$

$$+ \int_0^A \int_0^B \int_s^A \psi_1(\sigma)\{\Delta_{u^\sim} g_{12}(\sigma,B,s,t,y(s,t),u(s,t)) + [\nabla_y g_{12}(\sigma,B,s,t,y(s,t),u(s,t)) +$$

$$+ \nabla_y \Delta_{u^\sim} g_{12}(\sigma,B,s,t,y(s,t),u(s,t))]\Delta y(s,t) + \mathbf{o}_{y,1,12}(\| \Delta y(s,t) \|)\}\,d\sigma\,dt\,ds +$$

$$+ \int_0^A \int_0^B \psi_2(t)\{\Delta_{u^\sim} g_1(A,t,s,y(s,t),u(s,t)) +$$

$$+ [\nabla_y g_1(A,t,s,y(s,t),u(s,t)) + \nabla_y \Delta_{u^\sim} g_1(A,t,s,y(s,t),u(s,t))]\Delta y(s,t) +$$

$$\mathbf{o}_{y,2,1}(\| \Delta y(\sigma,t) \|)\}\,dt\,ds +$$

$$+ \int_0^B \int_t^B \psi_2(\tau)\{\nabla_y g_2(A,\tau,t,y(A,t),u^*(A,t))\Delta y(A,t) + \mathbf{o}_{y,2,2}(\| \Delta y(A,t) \|)\}\,d\tau\,dt +$$

$$+ \int_0^B \int_0^A \int_t^B \psi_2(\tau)\{\Delta_{u^\sim} g_{12}(A,\tau,s,t,y(s,t),u(s,t)) + [\nabla_y g_{12}(A,\tau,s,t,y(s,t),u(s,t)) +$$

$$+ \nabla_y \Delta_{u^\sim} g_{12}(A,\tau,s,t,y(s,t),u(s,t))]\Delta y(s,t) + \mathbf{o}_{y,2,12}(\| \Delta y(s,t) \|)\}\,d\tau\,ds\,dt +$$

$$+\int_0^A \psi_0\{\nabla_y g_1(A,B,s,y(s,B),u^*(s,B))\Delta y(s,B)+\mathbf{o}_{y,0,1}(\|\Delta y(s,B)\|)\}\,ds+$$

$$+\int_0^B \psi_0\{\nabla_y g_2(A,B,s,y(A,t),u^*(A,t))\Delta y(A,t)+\mathbf{o}_{y,0,2}(\|\Delta y(A,t)\|)\}\,dt+$$

$$+\int_0^A\int_0^B \psi_0\{\Delta_{\tilde{u}}\, g_{12}(A,B,s,t,y(s,t),u(s,t))+[\nabla_y g_{12}(A,B,s,t,y(s,t),u(s,t))+$$

$$+\nabla_y\Delta_{\tilde{u}}\, g_{12}(A,B,s,t,y(s,t),u(s,t))]\Delta y(s,t)+\mathbf{o}_{y,0,12}(\|\Delta y(s,t)\|)\}\,dt\,ds+$$

$$+\nabla_y G_0(y(A,B))\Delta y(A,B)+\mathbf{o}_{J,0,1}(\|\Delta y(A,B)\|)+$$

$$+\int_0^A \{\nabla_y G_1(s,y(s,B))\Delta y(s,B)+\mathbf{o}_{J,1,1}(\|\Delta y(s,B)\|)\}\,ds+$$

$$+\int_0^B \{\nabla_y G_2(t,y(A,t))\Delta y(A,t)+\mathbf{o}_{J,2,1}(\|\Delta y(A,t)\|)\}\,dt+$$

$$+\int_0^A\int_0^B \{\Delta_{\tilde{u}}\, G_{12}(s,t,y(s,t),u(s,t))+[\nabla_y G_{12}(s,t,y(s,t),u(s,t))+$$

$$+\nabla_y\Delta_{\tilde{u}}\, G_{12}(s,t,y(s,t),u(s,t))]\Delta y(s,t)+\mathbf{o}_{J,12,1}(\|\Delta y(s,t)\|)\}\,dt\,ds$$

--- (8.24)

Because of the Hamiltonian equations (8.3), some of the terms on the right-hand side of (8.24) make up E. This is seen from the following equations:

$$\int_0^A \int_0^B \nabla_y G_{12}(s,t,y(s,t),u(s,t))\Delta y(s,t)\,dt\,ds +$$

$$+\int_0^A \int_0^B [\psi_1(s)\nabla_y g_2(s,B,t,y(s,t),u(s,t)) + \psi_2(t)\nabla_y g_1(A,t,s,y(s,t),u(s,t)) +$$

$$+[\psi_1(s)\nabla_y g_2(s,B,t,y(s,t),u(s,t)) + \psi_2(t)\nabla_y g_1(A,t,s,y(s,t),u(s,t)) +$$

$$+\psi_0\nabla_y g_{12}(A,B,s,t,y(s,t),u(s,t))]\Delta y(s,t)\,dt\,ds +$$

$$+\int_0^A \int_0^B \int_s^A [\psi_{12}(\sigma,t)\nabla_y g_1(\sigma,t,s,y(s,t),u(s,t)) + \psi_1(\sigma)\nabla_y g_{12}(\sigma,B,s,t,y(s,t),u(s,t))]\cdot$$

$$\cdot\,\Delta y(s,t)\,d\sigma\,dt\,ds +$$

$$+\int_0^A \int_0^B \int_t^B [\psi_{12}(s,\tau)\nabla_y g_2(s,\tau,t,y(s,t),u(s,t)) + \psi_2(\tau)\nabla_y g_{12}(A,\tau,s,t,y(s,t),u(s,t))]\cdot$$

$$\cdot\,\Delta y(s,t)\,d\tau\,dt\,ds +$$

$$+\int_0^A \int_0^B \int_s^A \int_t^B \psi_{12}(\sigma,\tau)\nabla_y g_{12}(\sigma,\tau,s,t,y(s,t),u(s,t))\,\Delta y(s,t)\,d\tau\,d\sigma\,dt\,ds =$$

$$= \int_0^A \int_0^B \psi_{12}(s,t)\Delta y(s,t)\,dt\,ds$$

--- (8.25)

$$\int_0^A [\nabla_y G_1(s,y(s,B)) + \psi_0\nabla_y g_1(A,B,s,y(s,B),u(s,B))]\Delta y(s,B)\,ds +$$

$$+\int_0^A \int_s^A \psi_1(\sigma)\nabla_y g_1(\sigma,B,s,y(s,B),u(s,B))\Delta y(s,B)\,d\sigma\,ds = \int_0^A \psi_1(s)\Delta y(s,B)\,ds$$

--- (8.26)

$$\int_0^B [\nabla_y G_2(t,y(A,t)) + \psi_0\nabla_y g_2(A,B,t,y(A,t),u(A,t))]\Delta y(A,t)\,dt +$$

$$+\int_0^B \int_t^B \psi_2(\tau)\nabla_y g_2(A,\tau,t,y(A,t),u(A,t))\Delta y(A,t)\,d\tau\,dt = \int_0^B \psi_2(t)\Delta y(A,t)\,dt$$

--- (8.27)

$$\nabla_y G_0(y(A,B))\Delta y(A,B) = \psi_0 \Delta y(A,B)$$

--- (8.28)

The terms on the right-hand sides of (8.25) through (8.28) make up E. On the other hand, the terms on the right-hand side of (8.24), other than those that appear on the left-hand sides of (8.25) through (8.28), and other than those that are $\mathbf{o}(\| \Delta y \|)$, make up

$$\int_0^A \int_0^B \{\Delta_{u^\sim} \eta_{12}(s,t,y(s,t),\psi_0,\psi_1(.),\psi_2(.),\psi_{12}(.,.),u(s,t)) +$$
$$+ \nabla_y \Delta_{u^\sim} \eta_{12}(s,t,y(s,t),\psi_0,\psi_1(.),\psi_2(.),\psi_{12}(.,.),u(s,t))\Delta y(s,t)\} \, dt \, ds \, .$$

This proves (8.23), and thus also (8.11). This concludes Part I of the proof of the extremum principle for the control $u_{12}$.

Part II.

If $u_{12}^*$ is an optimal control, we denote by $y^*$ the corresponding state function, and by $(\psi_0^*, \psi_1^*, \psi_2^*, \psi_{12}^*)$ the corresponding solution of the Hamiltonian equations. If $(s_1, t_1) \in (0,A) \times (0,B)$ is a point of continuity of $u_{12}^*$, if $u^\wedge$ is an admissible value of the control $u_{12}$, we define a control $u^\sim_{\varepsilon\delta}$, for every $\varepsilon$ and $\delta$ positive and sufficiently small, by

$$u^\sim_{\varepsilon\delta}(s,t) = u^\wedge \ , \ \text{if } (s,t) \in (s_1, s_1+\delta) \times (t_1, t_1+\varepsilon);$$
$$u^\sim_{\varepsilon\delta}(s,t) = u_{12}^*(s,t) \ , \ \text{if } (s,t) \notin (s_1, s_1+\delta) \times (t_1, t_1+\varepsilon)$$

--- (8.29)

We set

$$D_1 := \{(s,t) : s_1 < s \le A,\ t_1 < t \le B\}\ ;$$
$$D^0_{\varepsilon\delta} := \{(s,t) : s_1 < s \le s_1 + \delta,\ t_1 < t \le t_1 + \varepsilon\}\ ;$$
$$D^1_{\varepsilon\delta} := \{(s,t) : s_1 < s \le s_1 + \delta,\ t_1 + \varepsilon < t \le B\}\ ;$$
$$D^2_{\varepsilon\delta} := \{(s,t) : s_1 + \delta < s \le A,\ t_1 < t \le t_1 + \varepsilon\}\ ;$$
$$D^3_{\varepsilon\delta} := \{(s,t) : s_1 + \delta < s \le A,\ t_1 + \varepsilon < t \le B\}$$

--- (8.30)

We denote by $y_{\tilde{\varepsilon\delta}}$ the solution of the state dynamics (8.1) with control $u_{\tilde{\varepsilon\delta}}$. Then $y_{\tilde{\varepsilon\delta}}$ solves the following problem:

$y_{\tilde{\varepsilon\delta}}(s,t) = y^*(s,t)$, for $(s,t) \notin D_1$;

$$y_{\tilde{\varepsilon\delta}}(s,t) = g_0(s,t,u(s,t)) + \int_0^{s_1} g_1(s,t,\sigma,y^*(\sigma,t),u^*(\sigma,t))d\sigma +$$

$$+\int_0^{t_1} g_2(s,t,\tau,y^*(s,\tau),u^*(s,\tau))d\tau + \int_0^{s_1}\int_0^{t_1} g_{12}(s,t,\sigma,\tau,y^*(\sigma,\tau),u^*(\sigma,\tau))d\tau d\sigma +$$

$$+\int_{s_1}^{s}\int_0^{t_1} g_{12}(s,t,\sigma,\tau,y^*(\sigma,\tau),u^*(\sigma,\tau))d\tau d\sigma + \int_0^{s_1}\int_{t_1}^{t} g_{12}(s,t,\sigma,\tau,y^*(\sigma,\tau),u^*(\sigma,\tau))d\tau d\sigma +$$

$$+\int_{s_1}^{s} g_1(s,t,\sigma,y_{\tilde{\varepsilon\delta}}(\sigma,t),u^{\wedge})d\sigma + \int_{t_1}^{t} g_2(s,t,\tau,y_{\tilde{\varepsilon\delta}}(s,\tau),u^{\wedge})d\tau +$$

$$+\int_{s_1}^{s}\int_{t_1}^{t} g_{12}(s,t,\sigma,\tau,y_{\tilde{\varepsilon\delta}}(\sigma,\tau),u^{\wedge})d\tau d\sigma, \text{ for } (s,t) \in D^0_{\varepsilon\delta};$$

$$y_{\tilde{\varepsilon\delta}}(s,t) = g_0(s,t,u^*(s,t)) + \int_0^{s_1} g_1(s,t,\sigma,y^*(\sigma,t),u^*(\sigma,t))d\sigma +$$

$$+\int_0^{t_1} g_2(s,t,\tau,y^*(s,\tau),u^*(s,\tau))d\tau + \int_0^{s}\int_0^{t_1} g_{12}(s,t,\sigma,\tau,y^*(\sigma,\tau),u^*(\sigma,\tau))d\tau d\sigma +$$

$$+\int_{s_1}^{s_1+\delta} g_1(s,t,\sigma,y_{\tilde{\varepsilon\delta}}(\sigma,t),u^{\wedge})d\sigma + \int_{t_1}^{t} g_2(s,t,\tau,y_{\tilde{\varepsilon\delta}}(s,\tau),u^*(s,\tau))d\tau +$$

$$+\int_{s_1}^{s_1+\delta}\int_{t_1}^{t} g_{12}(s,t,\sigma,\tau,y_{\tilde{\varepsilon\delta}}(\sigma,\tau),u^{\wedge})d\tau d\sigma + \int_{s_1+\delta}^{s} g_1(s,t,\sigma,y_{\tilde{\varepsilon\delta}}(\sigma,t),u^*(\sigma,t))d\sigma +$$

$$+\int_0^{s_1}\int_{t_1}^{t} g_{12}(s,t,\sigma,\tau,y^*(\sigma,\tau),u^*(\sigma,\tau))d\tau d\sigma + \int_{s_1+\delta}^{s}\int_{t_1}^{t} g_{12}(s,t,\sigma,\tau,y_{\tilde{\varepsilon\delta}}(\sigma,\tau),u^*(\sigma,\tau))d\tau d\sigma,$$

for $(s,t) \in D^1_{\varepsilon\delta}$;

$$y_{\tilde{\varepsilon\delta}}(s,t) = g_0(s,t,u(s,t)) + \int_0^{s_1} g_1(s,t,\sigma,y^*(\sigma,t),u^*(\sigma,t))\,d\sigma +$$

$$+ \int_0^{t_1} g_2(s,t,\tau,y^*(s,\tau),u^*(s,\tau))\,d\tau + \int_0^{s_1}\int_0^{t} g_{12}(s,t,\sigma,\tau,y^*(\sigma,\tau),u^*(\sigma,\tau))\,d\tau\,d\sigma +$$

$$+ \int_{s_1}^{s} g_1(s,t,\sigma,y_{\tilde{\varepsilon\delta}}(\sigma,t),u^*(\sigma,t))\,d\sigma + \int_{t_1}^{t_1+\varepsilon} g_2(s,t,\tau,y_{\tilde{\varepsilon\delta}}(s,\tau),\hat{u})\,d\tau +$$

$$+ \int_{t_1+\varepsilon}^{t} g_2(s,t,\tau,y_{\tilde{\varepsilon\delta}}(\sigma,t),u^*(s,\tau))\,d\tau + \int_{s_1}^{s}\int_0^{t_1} g_{12}(s,t,\sigma,\tau,y^*(\sigma,\tau),u^*(\sigma,\tau))\,d\tau\,d\sigma +$$

$$+ \int_{s_1}^{s_1+\delta}\int_{t_1}^{t_1+\varepsilon} g_{12}(s,t,\sigma,\tau,y_{\tilde{\varepsilon\delta}}(\sigma,\tau),\hat{u})\,d\tau\,d\sigma +$$

$$+ \int_{s_1}^{s}\int_{t_1+\varepsilon}^{t} g_{12}(s,t,\sigma,\tau,y_{\tilde{\varepsilon\delta}}(\sigma,t),u^*(\sigma,\tau))\,d\tau\,d\sigma\,, \text{ for } (s,t)\in D^2_{\varepsilon\delta}\,;$$

$$y_{\tilde{\varepsilon\delta}}(s,t) = g_0(s,t,u^*(s,t)) + \int_0^{s} g_1(s,t,\sigma,y^*(\sigma,t),u^*(\sigma,t))\,d\sigma +$$

$$+ \int_0^{t} g_2(s,t,\tau,y^*(s,\tau),u^*(s,\tau))\,d\tau + \int_0^{s_1}\int_0^{t_1} g_{12}(s,t,\sigma,\tau,y^*(\sigma,\tau),u^*(\sigma,\tau))\,d\tau\,d\sigma +$$

$$+ \int_0^{s_1}\int_{t_1}^{t} g_{12}(s,t,\sigma,\tau,y^*(\sigma,\tau),u^*(\sigma,\tau))\,d\tau\,d\sigma + \int_{s_1}^{s}\int_0^{t_1} g_{12}(s,t,\sigma,\tau,y^*(\sigma,\tau),u^*(\sigma,\tau))\,d\tau\,d\sigma +$$

$$+ \int_{s_1}^{s_1+\delta}\int_{t_1}^{t_1+\varepsilon} g_{12}(s,t,\sigma,\tau,y_{\tilde{\varepsilon\delta}}(\sigma,t),\hat{u})\,d\tau\,d\sigma + \int_{s_1}^{s_1+\delta}\int_{t_1+\varepsilon}^{t} g_{12}(s,t,\sigma,\tau,y_{\tilde{\varepsilon\delta}}(\sigma,t),u^*(\sigma,\tau))\,d\tau\,d\sigma +$$

$$+ \int_{s_1+\delta}^{s}\int_{t_1}^{t_1+\varepsilon} g_{12}(s,t,\sigma,\tau,y_{\tilde{\varepsilon\delta}}(\sigma,t),u^*(\sigma,\tau))\,d\tau\,d\sigma +$$

$$+ \int_{s_1+\delta}^{s}\int_{t_1+\varepsilon}^{t} g_{12}(s,t,\sigma,\tau,y_{\tilde{\varepsilon\delta}}(\sigma,t),u^*(\sigma,\tau))\,d\tau\,d\sigma\,, \text{ for } (s,t)\in D^3_{\varepsilon\delta}$$

--- (8.31)

Consequently, for the controls $u_{\tilde{\varepsilon\delta}}$ introduced above, the difference $\Delta y(s,t) := y_{\tilde{\varepsilon\delta}}(s,t) - y^*(s,t)$ can be non-zero only for $(s,t) \in D_1$. It is thus seen that $\Delta y(s,t)$ satisfies an equation

$$\Delta y(s,t) = \Delta_{u_{\tilde{\varepsilon\delta}}} g_0(s,t,u^*(s,t)) + \int_{s_1}^{s} \Delta g_1(s,t,\sigma,y^*(\sigma,t),u^*(\sigma,t))\,d\sigma +$$
$$+ \int_{t_1}^{t} \Delta g_2(s,t,\tau,y^*(s,\tau),u^*(s,\tau))\,d\tau + \int_{s_1}^{s}\int_{t_1}^{t} \Delta g_{12}(s,t,\sigma,\tau,y^*(\sigma,\tau),u^*(\sigma,\tau))\,d\tau\,d\sigma \quad \text{--- (8.32)}$$

where

$$\Delta g_1(s,t,\sigma,y^*(\sigma,t),u^*(\sigma,t)) := g_1(s,t,\sigma,y_{\tilde{\varepsilon\delta}}(\sigma,t),u_{\tilde{\varepsilon\delta}}(\sigma,t)) - g_1(s,t,\sigma,y^*(\sigma,t),u^*(\sigma,t)) \;;$$
$$\Delta g_2(s,t,\tau,y^*(s,\tau),u^*(s,\tau)) := g_2(s,t,\tau,y_{\tilde{\varepsilon\delta}}(s,\tau),u_{\tilde{\varepsilon\delta}}(s,\tau)) - g_2(s,t,\tau,y^*(s,\tau),u^*(s,\tau)) \;;$$
$$\Delta g_{12}(s,t,\sigma,\tau,y^*(\sigma,\tau),u^*(\sigma,\tau)) :=$$
$$= g_{12}(s,t,\sigma,\tau,y_{\tilde{\varepsilon\delta}}(\sigma,\tau),u_{\tilde{\varepsilon\delta}}(\sigma,\tau)) - g_{12}(s,t,\sigma,\tau,y^*(\sigma,\tau),u^*(\sigma,\tau)) \quad \text{--- (8.33)}$$

Each of the finite variations above can be expressed according to the equation

$$\Delta g_{12}(s,t,\sigma,\tau,y^*(\sigma,\tau),u^*(\sigma,\tau)) =$$
$$= \Delta_{y_{\tilde{\varepsilon\delta}}} g_{12}(s,t,\sigma,\tau,y^*(\sigma,\tau),u_{\tilde{\varepsilon\delta}}(\sigma,\tau)) + \Delta_{u_{\tilde{\varepsilon\delta}}} g_{12}(s,t,\sigma,\tau,y^*(\sigma,\tau),u^*(\sigma,\tau)) \quad \text{--- (8.34)}$$

which, of course, also applies to the functions $g_1$ and $g_2$, by replacing $\{(s,t,\sigma,\tau),(\sigma,\tau)\}$ by $\{(s,t,\sigma),(\sigma,t)\}$, $\{(s,t,\tau),(s,\tau)\}$, respectively. Also, we note that $\Delta_{u_{\tilde{\varepsilon\delta}}} g_{12}(s,t,\sigma,\tau,y^*(\sigma,\tau),u^*(\sigma,\tau))$ can be non-zero only on the set $D^0_{\varepsilon\delta}$.

If $B_1, B_2, B_{12}$ are Lipschitz constants for the functions $g_1, g_2, g_{12}$, respectively, with respect to the variable y, then we have

$$\| \Delta_{\tilde{y}_{\varepsilon\delta}} g_1(s,t,\sigma,y^*(\sigma,t),u_{\tilde{\varepsilon\delta}}(\sigma,t)) \| \le B_1 \| \Delta y(\sigma,t) \| ;$$

$$\| \Delta_{\tilde{y}_{\varepsilon\delta}} g_2(s,t,\tau,y^*(s,\tau),u_{\tilde{\varepsilon\delta}}(s,\tau)) \| \le B_2 \| \Delta y(s,\tau) \| ;$$

$$\| \Delta_{\tilde{y}_{\varepsilon\delta}} g_{12}(s,t,\sigma,\tau,y^*(\sigma,\tau),u_{\tilde{\varepsilon\delta}}(\sigma,\tau)) \| \le B_{12} \| \Delta y(\sigma,\tau) \| \qquad \text{--- (8.35)}$$

Consequently $\Delta y$ satisfies the integral inequality

$$\| \Delta y(s,t) \| \le a(s,t) + B_1 \int_{s_1}^{s} \| \Delta y(\sigma,t) \| d\sigma + B_2 \int_{t_1}^{t} \| \Delta y(s,\tau) \| d\tau +$$

$$+ B_{12} \int_{s_1}^{s} \int_{t_1}^{t} \| \Delta y(\sigma,\tau) \| d\tau d\sigma ;$$

$$a(s,t) := \Bigg[ \Delta_{\tilde{u}_{\varepsilon\delta}} g_0(s,t,u^*(s,t)) \chi((s,t) \in D^0_{\varepsilon\delta}) +$$

$$+ \int_{s_1}^{s \wedge (s_1+\delta)} \Delta_{\tilde{u}_{\varepsilon\delta}} g_1(s,t,\sigma,y^*(\sigma,t),u^*(\sigma,t)) d\sigma + \int_{t_1}^{t \wedge (t_1+\varepsilon)} \Delta_{\tilde{u}_{\varepsilon\delta}} g_2(s,t,\tau,y^*(s,\tau),u^*(s,\tau)) d\tau +$$

$$+ \int_{s_1}^{s \wedge (s_1+\delta)} \int_{t_1}^{t \wedge (t_1+\varepsilon)} \Delta_{\tilde{u}_{\varepsilon\delta}} g_{12}(s,t,\sigma,\tau,y^*(\sigma,\tau),u^*(\sigma,\tau)) d\tau d\sigma \Bigg] \qquad \text{--- (8.36)}$$

Here, $\chi$ is the truth function, i.e. if P is a logical statement, then $\chi(P) = 1$ if P is true, and $\chi(P) = 0$ if P is false. The symbol $\wedge$ denotes the operation of taking the minimum of two real numbers.
This concludes Part II.

Part III.

First we prove an inequality of Gronwall's type. Even though there exist, in the research literature, several extensions of the standard Gronwall inequality, we shall prove directly a result that we need for the present paper.
If $z(s,t), \varphi_0(s,t), \varphi_1(s,t,\sigma), \varphi_2(s,t,\tau), \varphi_{12}(s,t,\sigma,\tau)$ are piecewise continuous nonnegative functions and z satisfies

$$z(s,t) \le \varphi_0(s,t) + \int_{s_1}^{s} \varphi_1(s,t,\sigma) z(\sigma,t)\, d\sigma + \int_{t_1}^{t} \varphi_2(s,t,\tau) z(s,\tau)\, d\tau +$$

$$+ \int_{s_1}^{s} \int_{t_1}^{t} \varphi_{12}(s,t,\sigma,\tau) z(\sigma,\tau)\, d\tau\, d\sigma$$

--- (8.37)

and if $\zeta$ is the solution of

$$\zeta(s,t) = \varphi_0(s,t) + \int_{s_1}^{s} \varphi_1(s,t,\sigma) \zeta(\sigma,t)\, d\sigma + \int_{t_1}^{t} \varphi_2(s,t,\tau) \zeta(s,\tau)\, d\tau +$$

$$+ \int_{s_1}^{s} \int_{t_1}^{t} \varphi_{12}(s,t,\sigma,\tau) \zeta(\sigma,\tau)\, d\tau\, d\sigma$$

--- (8.38)

then

$$z(s,t) \le \zeta(s,t) \text{ on } [s_1, A] \times [t_1, B]$$

--- (8.39)

We shall use the convolution notation of section 5. (The fact that we now have integrations with starting point $(s_1, t_1)$, instead of (0, 0) that we had in section 5, does not affect the validity of the results.) We set $\varphi := (\varphi_1, \varphi_2, \varphi_{12})$, and we write (8.37) in the form

$$z \le \varphi_0 + \varphi \otimes_0 z$$

--- (8.40)

Because of the nonnegativity of all functions involved, we can iterate (8.40) to get

$$z \le \varphi_0 + \varphi \otimes_0 (\varphi_0 + \varphi \otimes_0 z) = \varphi_0 + \varphi \otimes_0 \varphi_0 + \varphi^{\otimes 2} \otimes_0 z$$

--- (8.41)

and, inductively,

$$z \le \varphi_0 + \sum_{k=1}^{n} \varphi^{\otimes k} \otimes_0 \varphi_0 + \varphi^{\otimes(n+1)} \otimes_0 z$$

--- (8.42)

The convergence analysis of section 5 shows that $\varphi^{\otimes(n+1)} \otimes_0 z \to 0$ as $n \to \infty$ in the uniform (sup) norm over $[s_1, A] \times [t_1, B]$. By passing to the limit in (8.42), we obtain

$$z \le \varphi_0 + \sum_{k=1}^{\infty} \varphi^{\otimes k} \otimes_0 \varphi_0$$

--- (8.43)

Now, the solution of (8.38) is

$$\zeta = \varphi_0 + \sum_{k=1}^{\infty} \varphi^{\otimes k} \otimes_0 \varphi_0$$

--- (8.44)

and consequently (8.43) is tantamount to the wanted inequality (8.40).

Next, we obtain more explicit information about the solution of (8.39) in the case in which all functions $\varphi_0(s,t), \varphi_1(s,t,\sigma), \varphi_2(s,t,\tau), \varphi_{12}(s,t,\sigma,\tau)$ are constants. Thus we consider the equation

$$\zeta(s,t) = A + B_1 \int_{s_1}^{s} \zeta(\sigma,t)\, d\sigma + B_2 \int_{t_1}^{t} \zeta(s,\tau)\, d\tau + B_{12} \int_{s_1}^{s} \int_{t_1}^{t} \zeta(\sigma,\tau)\, d\sigma\, d\tau$$

--- (8.41)

We seek a solution in the form

$$\zeta(s,t) = A \sum_{k=0}^{\infty} \sum_{\ell=0}^{\infty} C_{k\ell} (s - s_1)^k (t - t_1)^\ell$$

--- (8.42)

By substituting (8.42) into (8.41) and equating the coefficients of the same powers of $(s - s_1)$, $(t - t_1)$, we find

$$C_{k0}=\frac{B_1{}^k}{k!},\ C_{0\ell}=\frac{B_2{}^\ell}{\ell!},$$
$$C_{k+1,\ell+1}=\frac{B_1}{k+1}C_{k,\ell+1}+\frac{B_2}{\ell+1}C_{k+1,\ell}+\frac{B_{12}}{(k+1)(\ell+1)}C_{k\ell}$$

--- (8.43)

We set

$$B:=\max\left(|B_1|,|B_2|,\sqrt{\frac{|B_{12}|}{3}}\right)$$

--- (8.44)

We claim that

$$|\zeta(s,t)|\le|A|\exp\big(3B(s-s_1)(t-t_1)\big)$$

--- (8.45)

In order to prove (8.45), it suffices to show that

$$|C_{k\ell}|\le\frac{(3B)^{k+\ell}}{(k!)(\ell!)}$$

--- (8.46)

The inequality (8.46) is clearly true for $k=0$ or $\ell=0$. For other values of $(k,\ell)$, we use double induction on $(k,\ell)$, and we invoke (8.43):

$$|C_{k+1,\ell+1}|\le\frac{|B_1|}{k+1}|C_{k,\ell+1}|+\frac{|B_2|}{\ell+1}|C_{k+1,\ell}|+\frac{|B_{12}|}{(k+1)(\ell+1)}|C_{k\ell}|\le$$
$$\le\frac{|B_1|}{k+1}\frac{(3B)^{k+\ell+1}}{(k!)((\ell+1)!)}+\frac{|B_2|}{\ell+1}\frac{(3B)^{k+\ell+1}}{((k+1)!)(\ell!)}+\frac{|B_{12}|}{(k+1)(\ell+1)}\frac{(3B)^{k+\ell}}{(k!)(\ell!)}\le$$
$$\le 3B\frac{(3B)^{k+\ell+1}}{((k+1)!)((\ell+1)!)}=\frac{(3B)^{k+\ell+2}}{((k+1)!)((\ell+1)!)}$$

--- (8.47)

This proves (8.46) (and thus (8.45)), and (8.46) also proves the uniform convergence of the series (8.42).
This concludes Part III.

Part IV.

In the case in which $g_0, g_1, g_2$ are independent of u, the inequality (8.36) becomes

$$\| \Delta y(s,t) \| \le$$
$$\le \int_{s_1}^{s \wedge (s_1+\delta)} \int_{t_1}^{t \wedge (t_1+\varepsilon)} \Delta_{\tilde{u}_{\varepsilon\delta}} g_{12}(s,t,\sigma,\tau,y^*(\sigma,\tau),u^*(\sigma,\tau))\, d\tau\, d\sigma + B_{12} \int_{s_1}^{s} \int_{t_1}^{t} \| \Delta y(\sigma,\tau) \| d\tau\, d\sigma$$

--- (8.48)

If $C_0$ is a bound on $g_{12}(s,t,\sigma,\tau,y^*(\sigma,\tau),u)$, uniform with respect to $s,t,\sigma,\tau,u$, then the first double integral term in (8.48) is, in absolute value, $\le 2C_0\varepsilon\delta$. Consequently, by the results established in part III, we have

$$\| \Delta y(s,t) \| \le 2C_0\varepsilon\delta \exp(3M\varepsilon\delta) = \mathbf{O}_1(\varepsilon\delta)$$

--- (8.49)

abd similarly

$$\| \Delta y(A,t) \| = \mathbf{O}_2(\varepsilon\delta), \| \Delta y(s,B) \| = \mathbf{O}_3(\varepsilon\delta), \| \Delta y(A,B) \| = \mathbf{O}_4(\varepsilon\delta)$$

--- (8.50)

Now (8.11) becomes

$$\Delta J(y,u) = \iint_{D^0_{\varepsilon\delta}} \{\Delta_{\tilde{u}_{\varepsilon\delta}} \eta_{12}(s,t,y(s,t),\psi_0,\psi_1(.),\psi_2(.),\psi_{12}(.,.),u(s,t)) +$$
$$+ \nabla_y \Delta_{\tilde{u}_{\varepsilon\delta}} \eta_{12}(s,t,y(s,t),\psi_0,\psi_1(.),\psi_2(.),\psi_{12}(.,.),u(s,t)) \Delta y(s,t)\}\, dt\, ds +$$
$$+ \iint_{D_1} \mathbf{o}_{J,12}(\| \Delta y(s,t) \|)\}\, dt\, ds +$$
$$+ \int_{s_1}^{A} \mathbf{o}_{J,1}(\| \Delta y(s,B) \|)\, ds + \int_{t_1}^{B} \mathbf{o}_{J,2}(\| \Delta y(A,t) \|)\, dt + \mathbf{o}_{J,0}(\| \Delta y(A,B) \|)$$

--- (8.51)

In view of (8.49) and (8.50), we obtain from (8.51)

$$\Delta J(y,u) = \iint_{D^0_{\varepsilon\delta}} \{\Delta_{\tilde{u}_{\varepsilon\delta}} \eta_{12}(s,t,y(s,t),\psi_0,\psi_1(.),\psi_2(.),\psi_{12}(.,.),u(s,t))\,dt\,ds + \mathbf{o}(\varepsilon\delta)$$

--- (8.52)

For an optimal control $u^*$ and the corresponding state trajectory $y^*$, we have

$\Delta J(y^*,u^*) \geq 0$,

and, consequently,

$$\iint_{D^0_{\varepsilon\delta}} \{\Delta_{\tilde{u}_{\varepsilon\delta}} \eta_{12}(s,t,y^*(s,t),\psi_0{}^*,\psi_1{}^*(.),\psi_2{}^*(.),\psi_{12}{}^*(.,.),u^*(s,t))\,dt\,ds + \mathbf{o}(\varepsilon\delta) \geq 0$$

--- (8.53)

By taking $\delta = q\varepsilon$, for a fixed positive number q, we obtain

$$\iint_{D^0_{\varepsilon,q\varepsilon}} \{\Delta_{\tilde{u}_{\varepsilon,q\varepsilon}} \eta_{12}(s,t,y^*(s,t),\psi_0{}^*,\psi_1{}^*(.),\psi_2{}^*(.),\psi_{12}{}^*(.,.),u^*(s,t))\,dt\,ds + \mathbf{o}(\varepsilon^2) \geq 0$$

--- (8.54)

By invoking standard results on differentiability of integrals, we have, for almost all points $(s_1,t_1)$,

$$\iint_{D^0_{\varepsilon,q\varepsilon}} \{\Delta_{\tilde{u}_{\varepsilon,q\varepsilon}} \eta_{12}(s,t,y^*(s,t),\psi_0{}^*,\psi_1{}^*(.),\psi_2{}^*(.),\psi_{12}{}^*(.,.),u^*(s,t))\,dt\,ds =$$

$$= q\varepsilon^2 \Delta_{\tilde{u}_{\varepsilon,q\varepsilon}} \eta_{12}(s_1,t_1,y^*(s,t),\psi_0{}^*,\psi_1{}^*(.),\psi_2{}^*(.),\psi_{12}{}^*(.,.),u^*(s_1,t_1)) + \mathbf{o}_I(\varepsilon^2)$$

--- (8.55)

For points $(s_1,t_1)$ for which (8.55) holds, by dividing both sides of (8.54) by $q\varepsilon^2$ and taking the limit as $\varepsilon \to 0^+$, we obtain

$$\eta_{12}(s_1,t_1,y^*(s,t),\psi_0{}^*,\psi_1{}^*(.),\psi_2{}^*(.),\psi_{12}{}^*(.,.),u^*(s_1,t_1))-$$
$$-\eta_{12}(s_1,t_1,y^*(s,t),\psi_0{}^*,\psi_1{}^*(.),\psi_2{}^*(.),\psi_{12}{}^*(.,.),u_1)\geq 0$$

--- (8.55)

which is tantamount to the partial extremal principle for $u_{12}$ that we stated in section 7 above. ///